\documentclass[12pt]{amsart}
\usepackage{amscd,amssymb}
\newtheorem{thm}{Theorem}[section]
\newtheorem{corol}[thm]{Corollary} 
\newtheorem{lemma}[thm]{Lemma}
\newtheorem{prop}[thm]{Proposition}
\theoremstyle{definition}
\newtheorem{defin}[thm]{Definition}

\theoremstyle{remark}
\newtheorem{remark}[thm]{Remark}
\newtheorem{remarks}[thm]{Remarks}
\newtheorem{example}[thm]{Example}

\numberwithin{equation}{section}

\newcommand{\abs}[1]{\lvert#1\rvert}
\def\norm#1{\left\Vert#1\right\Vert}
\def\R {{\Bbb R}}

\def\I {{\Bbb I}}
\def\C {{\Bbb C}}
\def\N{{\Bbb N}}
\def\e{{\varepsilon}}
\def\Z {{\Bbb Z}}
\def\U{{\Bbb U}}
\def\F{{\mathrm{E}}}
\def\Un{{\mathcal{U}}}
\def\Iso{{\mathrm{Iso}}\,}

\def\St{{\mathrm St}\,}
\def\Homeo{{\mathrm{Homeo}}\,}

\def\diam{{\mathrm{diam}}}

\def\H{{\mathcal H}}
\def\s{{\mathbb S}}

\def\colon{{{:}\;}}
\def\endrem{}

\begin{document}

\title[Ramsey--Milman phenomenon]{Ramsey--Milman phenomenon, 
Urysohn metric spaces, and extremely amenable groups}

%    Information for first author
\author[V. Pestov]{Vladimir Pestov}
%    Address of record for the research reported here
\address{School of Mathematical and Computing Sciences,
Victoria University of Wellington, P.O. Box 600, Wellington,
New Zealand}
%    Current address
%\curraddr{}
\email{vova@mcs.vuw.ac.nz}
\urladdr{http://www.mcs.vuw.ac.nz/$^\sim$vova}
%    \thanks will become a 1st page footnote.

\thanks{{\it 2000 Mathematical Subject Classification.} 
Primary: 22F05. Secondary: 05C55, 22A05, 28C10, 43A05, 43A07, 51F99, 54H25.}

%\thanks{Preprint VUW-SMCS 00-8, revised.}
%\date{January 15, 2001}

\keywords{}
\begin{abstract} 
In this paper we further study links between
concentration of measure in topological transformation groups, 
existence of fixed points, and Ramsey-type theorems for metric spaces. 
We prove that whenever the group $\Iso(\U)$
of isometries of Urysohn's universal 
complete separable metric space $\mathbb U$, equipped with the
compact-open topology, acts upon an arbitrary compact space,
it has a fixed point.
The same is true if $\U$ is replaced with any generalized
Urysohn metric space $U$ that is sufficiently homogeneous.
Modulo a recent theorem by Uspenskij that
every topological group embeds into a topological
group of the form $\Iso(U)$, our result implies
that every topological group embeds into an extremely amenable group 
(one admitting an invariant multiplicative
mean on bounded right uniformly continuous functions). 
By way of the
proof, we show that every topological group is approximated by
finite groups in a certain weak sense. 
Our technique also results
in a new proof of the extreme amenability
(fixed point on compacta property) for infinite orthogonal groups. 
Going in the opposite direction,
we deduce some Ramsey-type theorems for metric subspaces of Hilbert
spaces and for spherical metric spaces from existing results
on extreme amenability of infinite unitary groups and groups
of isometries of Hilbert spaces.
\end{abstract}

\maketitle

%\setcounter{tocdepth}{1}
%\tableofcontents

\section{Introduction} The concept of amenability extends
from locally compact groups to arbitrary topological groups, and 
an interesting observation of recent times is that under such a
transition the concept `gains in strength' in that a number of concrete
infinite-dimensional groups of importance satisfy 
a reinforced version of amenability such as locally compact groups 
cannot possibly have. 

Definitions of amenability equivalent in the locally compact 
case diverge already for some of the most common
infinite-dimensional topological groups \cite{dlH}.
Nevertheless, the following choice has become standard 
\cite{Pat,Aus}: 
call a topological group $G$ {\bf amenable} if every continuous
affine action of $G$ on a convex compact set has a fixed point.
Equivalently, there is a left invariant 
mean on the space $C^b_\Rsh(G)$
of all bounded right uniformly continuous functions on $G$.
This concept is in particular given substance by the following
result due to de la Harpe \cite{dlH2}:
a von Neumann algebra $A$ is injective if and only if the unitary group
$U(A)$ equipped with the ultraweak topology is amenable. 
(Cf. also \cite{Pat2}.) 
Such results suggest that 
namely the above definition and not, for example, the one
calling for an invariant mean on {\it all} bounded continuous
functions on $G$, is the `proper' choice.

In particular, a topological group $G$ is amenable if it has 
a fixed point in {\it every} compact space it acts upon.
Such topological groups are said to have the 
{\bf fixed point on compacta property}
({\bf f.p.c.}) \cite{Gl}, or else 
called {\bf extremely amenable}, in the spirit of \cite{Gra}
where the concept was applied to discrete semigroups.
The condition is equivalent to the existence of a left 
invariant {\bf multiplicative} mean on $C^b_\Rsh(G)$. 

At the first sight, the latter 
property seems to be far too restrictive
to be observed {\it en masse.} 
In particular, according to a well-known
theorem of Veech \cite{V}, no locally compact group has the fixed point
on compacta property. 
(For discrete groups, this was previously noted in \cite{Ell}.) 
Historically the first examples of extremely amenable groups \cite{HC, Ba},
difficult to construct, looked like genuine pathologies. 

Nevertheless, in recent times it was shown that a number of
well-known `massive' topological groups possess the fixed
point on compacta property, among them
\begin{list}{$\bullet$}{}
\item
the unitary group $U(\H)$ (and the orthogonal group
$\mathrm{O}(\H)$) of an infinite-dimensional
Hilbert space with the strong operator topology
(Gromov and Milman \cite{GrM}), 
\item
the group $L_1(X,U(1))$ of measurable maps from a non-atomic
Lebesgue space to the circle group, equipped with 
the $L_1$-metric (Glasner \cite{Gl} and independently, unpublished,
Furstenberg and B. Weiss),
\item groups $\Homeo_+({\mathbb I})$ and 
$\Homeo_+({\mathbb R})$ of orientation-preserving homeomorphisms
with the compact-open topology (the present author \cite{P1}),
\item groups of
measure-preserving automorphisms of standard sigma-finite measure
spaces with the strong topology 
(Giordano and the present author \cite{GP}).
\end{list}

The technique used to establish the fixed point on compacta property
in the above examples 
has been either that of concentration of measure 
on high-dimensional structures
(pioneered in this context by Gromov and Milman \cite{GrM}), or
else infinite Ramsey theory, as in \cite{P1}.

In this article we isolate a new and  
vast class of topological groups
with the fixed point on compacta property: they are
groups of isometries of very regular and highly homogeneous objects,
the (generalized) Urysohn metric spaces.

Universal metric spaces were introduced by
Urysohn in the 20's \cite{Ur1,Ur} 
and investigated mostly in the separable case.
In particular, there is, up to an isometry, only one complete separable
Urysohn metric space, which we will denote by $\U$.
For a long time Urysohn spaces remained little known outside of general
topology, and the most important advances at that period were due to 
Kat\u etov \cite{Kat}, who had made the structure of the space $\U$ more
transparent, and Uspenskij \cite{U1}, who had proved 
that the group of isometries 
$\Iso(\U)$ with the compact-open topology 
forms a universal second-countable topological group.
Uspenskij's construction was later used by Gao and Kechris 
\cite{GK} to deduce,
among others, the following result: every Polish topological group
is the group of {\it all} isometries of a suitable separable
complete metric space. 
Recently the
Urysohn spaces were linked to wider issues in geometry and analysis,
particularly by Vershik who has for example 
shown \cite{Ver} that the completion of 
the set of integers equipped with a
`sufficiently random' metric is almost surely isometric to $\U$.
A further discussion of the space $\U$ and its links with geometry
is to be found in Gromov's book \cite{Gr}.

We shall prove that the group $\Iso(\U)$ has the fixed point on
compacta property (Theorem \ref{fpc}), 
and moreover the same is true of isometry groups $\Iso(U)$
of all sufficiently homogeneous generalized
(non-separable) Urysohn spaces $U$ (Theorem \ref{generurysohn}). 
According to a recent result by Uspenskij \cite{U2}, 
every topological group
is contained, as a subgroup, in the group of isometries of such a
generalized Urysohn space. 
The two results combined imply that 
extreme amenability is, in a sense, ubiquitous:
every topological group embeds,
as a topological subgroup, into a topological group with the
fixed point on compacta property (Corollary \ref{main}).

It is known since the work of de la Harpe \cite{dlH}
that a closed subgroup of an amenable topological group need not be
amenable, unlike in the locally compact case.
The reported results take this observation to its extreme.
The possibility of such a development 
was conjectured in our paper \cite{P1}. 

The proof of extreme amenability of the group $\Iso(\U)$
applies the technique of concentration of measure, and
by way of proof we establish the following generalization of
a result due to Glasner and Furstenberg--Weiss: the group of
all measurable maps from a non-atomic Lebesgue measure space
to an amenable locally compact group $G$, equipped with the topology
of convergence in measure (known as
the Hartman--Mycielski extension of $G$, \cite{HM}),
has the fixed point on compacta property (Theorem \ref{gener}). 
Another component
of the proof is the following, apparently new, result
(Theorem \ref{approx}): every
group of isometries of a metric space can be approximated in a
certain weak sense with finite groups of isometries of suitable
metric spaces. 
In the second-countable case 
the result can be interpreted as a statement on
approximation of {\it topological} groups:
every Polish group is the limit of a net of finite groups
in the space of all closed 
subgroups of the group $\Iso(\U)$ (Corollary \ref{polish}).

Our methods lead to
a new proof of the fixed point on compacta property for the infinite
orthogonal groups with the strong topology, which does not use
advanced geometric tools such as Gromov's isoperimetric inequality.
(Subsection \ref{newp}.)

In order to extend the result
on extreme amenability to the groups of isometries $\Iso(U)$ of
generalized Urysohn metric spaces $U$, we recast the fixed
point on compacta property of the full isometry group of a sufficiently
homogeneous metric space $X$ as a Ramsey-type result for the space
$X$ itself (Theorem \ref{ramsey}). 
As a corollary, if two metric spaces,
$X$ and $Y$, are both $\omega$-homogeneous and have, up to
isometry, the same finite metric subspaces, then the groups
$\Iso(X)$ and $\Iso(Y)$ have the fixed point on
compacta property (or otherwise) simultaneously 
(Theorem \ref{ifandonly}). 

As another application of this technique, 
we show that the groups of isometries
of the universal discrete metric spaces \cite{Gr} do not have the fixed
point on compacta property (Theorem \ref{doesnot}). 

The equivalence between the fixed point on compacta property
of isometry groups and Ramsey-type results for metric spaces can be
exploited in the other direction as well, and thus we
deduce some `approximate' Ramsey-type
results for both spherical and Euclidean metric spaces 
(Subsection \ref{last}).

\subsection*{Acknowledgements}
The investigation grew out of stimulating
discussions with Mikha\"\i l Gromov and Vladimir Uspenskij, who have
both independently conjectured extreme amenability of 
the isometry groups of Urysohn spaces.
Moreover, Theorems \ref{approx} and \ref{rdm} are rigorous incarnations of
two further conjectures made by VU.
My thanks go to both named mathematicians,
to MG also for hospitality extended at IH\'ES in September 1999.
I am grateful to Alekos Kechris and to Eli Glasner for spotting
badly flawed fragments in the earlier versions of the article, and to
Michael Megrelishvili for a useful remark.
I am thankful to Sina Greenwood for the
invitation to give the opening talk at the 2000 Devonport Topology
Festival, which
had rekindled my work on the topic, and to Michael Cowling
for his hospitality at the Department of Pure Mathematics of the
University of New South Wales in April 2000, 
where part of the work had been done.
The research was partially supported by 
the Royal Society of New Zealand through the Marsden Fund
grant VUW703, by the Australian Research Council through the
large research grant `Random algebraic structures,' and by
the Institut des Hautes \'Etudes Scientifiques through its visitors programme.

\section{Concentration of measure in Hartman--Mycielski groups}
\subsection{}
Our starting point is the following result, mentioned in the
Introduction.

\begin{thm}[Glasner \cite{Gl}; Furstenberg--B. Weiss, unpublished]\hfill
\label{gfw}
The group $L_1(X,\mathrm{U}(1))$ 
of all measurable maps from a
nonatomic Lebesgue space to the circle rotation group, equipped
with the $L_1$-metric, has the fixed point on compacta property.
\end{thm} 

On two occasions in this article, including the proof of one
of our main theorems, we will invoke 
suitable modifications of the above result, and
it seems appropriate to state a far-reaching generalization of 
Theorem \ref{gfw}, even if we shall never use its full power. 
\par 
In the above form the result does not extend too far: suffice
to consider the additive group of the Banach space $L_1(X)=L_1(X,\R)$, with its
wealth of continuous characters. 
However, it is not a particular metric
on the group but rather the topology it generates that matters, and 
the topology generated by the $L_1$-metric on the group
$L(X,\mathbb{T})$ is that of convergence in measure. 
(This is true of
every $L_p$-metric, $1\leq p <\infty$, on the same group.) 
This observation leads us to state the following generalization 
of Glasner--Furstenberg--Weiss theorem.

\begin{thm}
Let $G$ be an amenable locally compact group and let $X$ be a
non-atomic Lebesgue measure space. 
Then the group $L_0(X,G)$ of all measurable maps from $X$ to $G$, 
equipped with the topology of convergence in measure,
has the fixed point on compacta property 
{\rm (}is extremely amenable{\rm )}.
\label{gener}
\end{thm}

\begin{remark}
The topological groups of the form $L_0(X,G)$, where the subscript
`$0$' stands for the topology of convergence in measure, 
had apparently been first considered by Hartman and Mycielski \cite{HM},
who had observed that $L_0(X,G)$ contains
$G$ as a topological subgroup (formed by all constant functions)
and is path-connected and locally path--connected. 
Later it was shown by Keesling \cite{Kee}
that if $G$ is separable metrizable, then the
Hartman--Mycielski extension $L_0(X,G)$ is homeomorphic to the
separable Hilbert space.
The correspondence $G\mapsto L_0(X,G)$ 
determines a (covariant) functor from the category of all topological
groups and continuous homomorphisms to itself, and Theorem
\ref{gener} says that the Hartman--Mycielski functor transforms
amenable locally compact groups into extremely amenable
topological groups.
\endrem\end{remark}

The following particular case (where $G=\R$ or $\C$)
seems to be of interest.

\begin{corol}
The {\rm [}underlying topological group of {\rm ]} 
the topological vector
space $L_0(X)$ of all measurable functions
on a non-atomic Lebesgue measure space $X$,
equipped with the topology of convergence in measure, has the
fixed point on compacta property. \qed
\end{corol}

\begin{remark}
The above result is similar to the one from \cite{HC} where the
space $L_0(X)$ was equipped with the topology of convergence
in a suitably chosen, the so-called
pathological submeasure (a subadditive set function) on $X$.
As a result, the abelian topological group from \cite{HC} has an even
stronger property than just extreme amenability: it admits no
strongly continuous unitary representations.
Notice that each of the groups of the form $L_0(X,G)$ from Theorem
\ref{gener} admits a faithful strongly continuous unitary representation
in the Hilbert space $L_2(X,L_2(G))$. (This extends an observation made
in \cite{Gl} for $G=U(1)$.)
\endrem\end{remark}

Our proof of Theorem \ref{gener} relies, similarly to that of Theorem
\ref{gfw}, on the technique
of concentration of measure on high-dimensional structures.
However, the concept of a L\'evy group \cite{GrM,Gl} becomes too
narrow and has to be somewhat extended. 
We believe that
this extension goes sufficiently far to be of interest on its own.
(Though we find it useful, to replace metrics with 
uniform structures, this is not what our generalization is about.)

\subsection{}
If $X=(X,\Un_X)$ is a uniform space, then the uniform 
(induced) topology on $X$
gives rise to a Borel structure and thus one can speak of Borel
measures on $X$. 
The following is a straightforward
adaptation of the by now classical concept \cite{GrM,MS,M2,Ta,Gr,Gl}.

\begin{defin}
Let $(\mu_\alpha)$ be a net of probability measures on 
a uniform space $(X,\Un_X)$.
Say that the net $(\mu_\alpha)$
has the {\bf L\'evy concentration property,} or simply
{\bf concentrates} (in $X$), if
whenever $A_\alpha\subseteq X$ are Borel subsets with the
property 
\[
\liminf_\alpha\mu_\alpha(A_\alpha)>0,
\] 
one has for every entourage of the diagonal $V\in\Un_X$
\[
\mu_\alpha(V[A_\alpha])\to 1.
\]
\endrem
\end{defin}
(Here, as usual, $V[A]=\{x\in X \mid \exists a\in A, ~(x,a)\in V\}$ denotes the
$V$-neighbourhood of $A$.)

\begin{lemma}
Let $f\colon X\to Y$ be a uniformly continuous map between two
uniform spaces, and let $(\mu_\alpha)$ be a net of Borel measures
on $X$. 
If $(\mu_\alpha)$ concentrates, then 
the net $(f_\ast(\mu_\alpha))$ of
push-forward measures on $Y$ concentrates as well. \qed
\label{conc-trans}
\end{lemma}

Let $G$ be a group of uniform isomorphisms of a uniform space $X$.
A compactification $K$ of $X$ is called {\bf uniform}
if the corresponding mapping $i\colon X\to K$ is uniformly continuous,
and {\bf equivariant} (in full,
$G${\bf -equivariant}) if $G$ acts on $K$ by homeomorphisms in
such a way that $i$ commutes with the action.
The maximal
uniform compactification of a uniform space $X$, known as the
{\bf Samuel compactification} of $X$ and which we denote by $\sigma X$,
is the Gelfand space of the commutative $C^\ast$-algebra formed by
all bounded uniformly continuous complex-valued functions on $X$.
The Samuel compactification $\sigma X$
is equivariant no matter what the acting group $G$ is, because every
uniform homeomorphism $X\to X$ extends to a self-homeomorphism
$\sigma X\to\sigma X$ due to universality.

It is convenient to state explicitely the following result, which is
in essence folk's knowledge in theory of topological
transformation groups.
(Cf. \cite{M2,V} etc.)

\begin{thm}
\label{folk}
Let $G$ be a group of uniform isomorphisms of a uniform space $X$.
The following conditions are equivalent.
\begin{enumerate}
\item[(i)] Every $G$-equivariant uniform compactification of $X$ has
a fixed point.
\item[(ii)] For every bounded uniformly continuous function $f$ from $X$
to a finite-dimensional Euclidean space, every $\e>0$ and 
every finite collection
$g_1,g_2,\dots,g_n\in G$, there is an $x\in X$ with
$\abs{f(x)-f(g_ix)}<\e$ for all $i=1,2,\dots,n$.
\item[(iii)] 
For every finite cover $\gamma$ of $X$, every $V\in\Un_X$ and 
every finite collection $g_1,g_2,\dots,g_n\in G$, there is
an $A\in\gamma$ such that 
\[
\cap_{i=1}^n V[g_iA]\neq\emptyset.
\]
\end{enumerate}
\end{thm}

\begin{proof} 
(i) $\Rightarrow$ (ii): notice that every bounded uniformly continuous
function extends over the Samuel compactification. \par
(ii) $\Rightarrow$ (iii): choose a bounded uniformly
continuous pseudometric $d$ on $X$ subordinated to $V$
(that is, $d(x,y)<1\Rightarrow (x,y)\in V$) and 
apply the condition (ii) to the function
from $X$ to $\R^{\abs\gamma}$ whose components are distance functions
$x\mapsto d(x,A)$, $A\in\gamma$, with $\e=1$. \par
(iii) $\Rightarrow$ (i): see \cite{P4}, Proposition 2.1
(which is, in its turn,
an adaptation of an argument from Section 4 in \cite{M1}.)
\end{proof}

\begin{remark}
At this point we do not concern ourselves with a 
topology on the acting group $G$, and it may well happen that
if $G$ is a topological group acting on the uniform space $X$
continuously, the extension of the action
to an equivariant compactification of $X$ is discontinuous.
As an example, consider as $G$ the unitary group $U(l_2)$ with
the strong topology,
and as $X$ the unit sphere $\s^\infty$ in $l_2$ with the metric
uniformity.
The action of $U(l_2)$ on the Samuel compactification
of the sphere is continuous if $U(l_2)$ is equipped with
the uniform operator topology, but not the strong one.
\endrem
\end{remark}

A subset $B$ of a uniform space $X$ is called
{\bf uniformly open} if $B=V[A]$ for some $A\subseteq X$ and
$V\in\Un_X$. 

\begin{defin}
Let us say that two nets of probability measures,
$(\mu_\alpha)$ and $(\nu_\alpha)$, 
on the same uniform space $X$ are {\bf asymptotically proximal} 
if for every uniformly open subset $B$ one has
\[
\limsup_\alpha\abs{\mu_\alpha(B)-\nu_\alpha(B)}<1
\]
\endrem\end{defin}

\begin{remark}
Two nets as above will in particular be asymptotically proximal if
$\liminf_\alpha(\mu_\alpha\wedge\nu_\alpha)(X)>0$.
For instance, this is so if the restrictions of 
$\mu_\alpha$ and $\nu_\alpha$ coincide on some Borel subsets
$(A_\alpha)$, whose measures are uniformly (in $\alpha$)
bounded away from zero.
\label{instance}
\endrem\end{remark}

The following apparently subsumes
all the previously known results of the type 
({\it concentration of measure})
$\Rightarrow$ ({\it existence of a fixed point})
\cite{GrM,M1,M2,Gl,P4}.

\begin{thm}
Let a group $G$ act on a uniform space 
$X=(X,\Un)$ by uniform isomorphisms.
Suppose there is a net 
$(\mu_\alpha)$ of probability measures on $X$ 
such that 
\begin{list}{-}{}
\item $(\mu_\alpha)$ concentrates in $X$,
\item for every $g\in G$ the nets $(\mu_\alpha)$ and
$(g\ast\mu_\alpha)$ are asymptotically proximal.
\end{list}
Then every equivariant uniform compactification of the 
$G$-space $X$ has a fixed point.
\end{thm}

\begin{remark}
The second condition is a rather weak 
invariance-type property for a family of measures, and its advantage is
being easier to verify.
If we require all measures $\mu_\alpha$ to be {\bf eventually invariant}
(that is, for every $g\in G$ one has $\mu_\alpha=g\ast\mu_\alpha$
for sufficiently large $\alpha$) and
compactly-supported, then we recover the concept of a L\'evy
transformation group from \cite{M1}.
The above stated theorem 
allows for a unified approach to a number of previously known
results, such as a link between amenability of unitary representations 
and the concentration property of unit spheres \cite{P4},
which we will not be addressing here.
\endrem\end{remark}

\begin{proof} 
Let $\gamma$ be a finite cover of $X$, let $V\in\Un_X$, 
and let $g_1,\dots,g_n\in G$ be arbitrary. 
Find an entourage of the diagonal
$W\in\Un_X$ with $W\circ W\subseteq V$.
At least one element of $\gamma$, denote it
by $A$, satisfies the property 
\[
\limsup_\alpha\mu_\alpha(A)\geq \abs\gamma^{-1}.
\]
By proceeding to a subnet if necessary, we may assume
without loss in generality that
\[
\liminf_\alpha\mu_\alpha(A)\geq \abs\gamma^{-1}.
\]
In view of the assumed 
concentration property of the measures $(\mu_\alpha)$,
\[
\lim_\alpha\mu_\alpha(W[A])=1,
\]
and by force of the second assumption, one has for every $i$
\[
\liminf_\alpha (g_i\ast\mu_\alpha)(W[A])>0.
\]
By Lemma \ref{conc-trans}, each of the nets of measures
$(g_i\ast\mu_\alpha)$, $i=1,2,\dots,n$ concentrates, and
consequently 
\[
\lim_\alpha (g_i\ast\mu_\alpha)(W[W[A]])=1.
\]
Since $W\circ W\subseteq V$, one has
\[
\lim_\alpha \mu_\alpha(g_i V[A])=1.
\]
It is therefore possible to choose an $\alpha$ so large 
that each of the numbers $\mu_\alpha(g_1 V[A])$, $\dots$,
$\mu_\alpha(g_n V[A])$ is greater than $1-\frac{1}{n}$. 
It follows that the intersection of all the translates of $V[A]$
by elements $g_i$, $i=1,2,\dots,n$ is non-empty, and application
of Theorem \ref{folk} finishes the proof.
\end{proof}

\subsection{}
Recall that the {\bf right uniform structure} of a topological
group, $\Un_\Rsh(G)$, has as a basis the entourages of diagonal 
of the form
\[
V_\Rsh=\{(x,y)\in G\times G \mid xy^{-1}\in V\},
\]
where $V$ runs over a neighbourhood basis of $e$ in $G$. 
The Samuel compactification of the right uniform space
$G_\Rsh=(G,\Un_\Rsh(G))$ is a compact $G$-space, known as the
{\bf greatest ambit} of $G$ and denoted by $\mathcal{S}(G)$.
(Cf. \cite{Te,Br,Aus,P2}.)
The greatest ambit possesses a distinguished
point (the image of identity of $G$, which we will still denote
$e$), whose orbit is everywhere
dense in it.
This object has the following universal property:
whenever $X$ is a compact $G$-space and $x_0\in X$, there is
a unique morphism of $G$-spaces from $\mathcal{S}(G)$ to $X$ taking
$e$ to $x_0$.
It follows that 
a topological group $G$ has the fixed point on compacta property
if and only if there is a fixed point in the greatest ambit
$\mathcal{S}(G)$.

\begin{corol}
Let $G$ be a topological group.
Suppose there is a net 
$(\mu_\alpha)$ of probability measures on $G$ such that,
with respect to the right uniform structure $\Un_\Rsh(G)$,
\begin{list}{-}{}
\item $(\mu_\alpha)$ concentrates,
\item for every $g\in G$ the nets $(\mu_\alpha)$ and
$(g\ast\mu_\alpha)$ are asymptotically proximal.
\end{list}
Then $G$ has the fixed point on compacta property.
\label{levy}
\qed
\end{corol}

\begin{remark}
A topological group $G$ is called a {\bf L\'evy
group} if it contains a family of compact 
subgroups, directed by inclusion and having everywhere dense union, 
such that the corresponding normalized 
Haar measures, $\mu_\alpha$, concentrate in $G_\Rsh$. 
This concept was used
as means to deduce the existence of fixed points for group actions
on compacta by Gromov and Milman \cite{GrM}; see also \cite{Gl, P2}. 
L\'evy groups satisfy a stronger property than the 
the second assumption of 
Corollary \ref{levy}: the measures $\mu_\alpha$ are eventually invariant. 
\endrem\end{remark}

\subsection{}
Let $X=(X,{\mathcal U}_X)$ be a uniform space. 
Denote by $L(\I,X)$ the collection of all 
Borel-measurable maps $f\colon \I\to X$ equipped with the
uniform structure of convergence in measure.
The standard basic entourages of diagonal are of the form
\begin{eqnarray*}
[V,\e] &:=&\{ (f,g)\in L(\I,X)\times L(\I,X) \colon \\
&& \mu\{x\in\I\colon (f(x),g(x))\notin V\} <\e\},
\end{eqnarray*}
where $V\in {\mathcal U}_X$ and $\e>0$.
This uniformity induces a topology on $L(\I,X)$, whose standard basic
neighbourhoods of a given function 
$f\colon \I\to X$ are 
\[
[V,\e,f] :=\{ g\in L(\I,X) \colon \mu\{x\in\I\colon (f(x),g(x))\notin V\}
<\e\},
\]
where $V\in {\mathcal U}_X$ and $\e>0$.
(Notice that the knowledge
of topology on $X$ alone does not suffice: to talk of convergence
in measure, it is necessary to have a uniform structure, 
for instance, one defined by a metric on $X$, or else the unique
compatible uniform structure in case $X$ is compact.)

If $G$ is a Hausdorff topological group, then so is $L(\I,G)$.
In this case, the standard neighbourhoods of identity are of the form 
\[
[V,\e] :=\{ g\in L(\I,X) \colon \mu\{x\in\I\colon g(x)\notin V\}<\e\},
\]
where $V$ is a neighbourhood of $e_G$ in $G$ and $\e>0$.

Now suppose that $X=(X,\rho)$ is a metric space.
Let us agree on the canonical
choice of the metric generating the uniformity of convergence in
measure (and the corresponding topology) 
on $L_0(\I,X)$, as follows: 
if $\lambda>0$ is an arbitrary (but fixed) number, then set
\begin{equation}
\mathrm{me}_\lambda(f,g)
=\inf\{ \e>0 \mid \mu\{x\in\I\colon \rho(f(x),g(x))
>\e\}<\lambda\e\}.
\label{metric}
\end{equation}
Such metrics for different $\lambda>0$ are all equivalent.
(Cf. \cite{Gr}, p. 115.)

\begin{defin}
An action of a topological group $G$ on a
uniform space $X=(X,\Un_X)$ by uniform isomorphisms is called {\bf bounded} 
\cite{dV} (or {\bf motion equicontinuous} \cite{GH})
if for every entourage $U\in\Un_X$ one can find
a neighbourhood $W\ni e_G$ such that for every $x\in X$, 
\[
W\cdot x\subseteq U[x].
\]
\label{defuca}
\endrem\end{defin}

\begin{remarks}
\label{uca}
1. Every bounded action is continuous. 
[If $g\in G$, $x\in X$, and a neighbourhood ${\mathcal O}\ni g\cdot x$
are arbitrary, select an $U\in \Un_X$ with 
$(U\circ U)[x]\subseteq {\mathcal O}$ and a neighbourhood 
$W\ni e_G$ with $W\cdot y\subseteq U[y]$ for all $y\in X$.
Then $W\cdot U[x] \subseteq U^2[x]\subseteq {\mathcal O}$.]
\par
2. The converse is not true.
[For example, 
the standard action of the unitary group $U(l_2)$ with the 
{\it strong} operator topology on the unit sphere $\s^\infty$
equipped with the metric uniformity is continuous, but not
bounded. This action becomes bounded if $U(l_2)$
is equipped with the {\it uniform} operator topology.]
\par
3. However, a continuous action of a topological 
group $G$ on a {\it compact} space $X$ (equipped with the unique
compatible uniformity) is always bounded.
This fact is well-known (and easily verified).
\par
4. The left action of a topological group $G$ on the {\it right}
uniform space $G_\Rsh$ is bounded, but in general the same is not
true of the left action of $G$ on the {\it left} uniform space of $G$.
\endrem
\end{remarks}

\begin{lemma}
\label{L}
If a topological group $G$ acts by isometries
on a metric space $X$, then the topological
group $L_0(\I,G)$ acts by isometries
on the metric space $L(\I,X)$ equipped with the metric (\ref{metric}),
where the action is defined pointwise:
\[
(g\cdot f)(x):=g(x)\cdot f(x), ~g\in L(\I,G), ~f\in L(\I,X).
\] 
If in addition the action of $G$ on $X$ is bounded
{\rm (}for example $X$ is compact{\rm )},
then the action of $L_0(\I,G)$ on $L_0(\I,X)$
is continuous.
\end{lemma}

\begin{proof} The first statement is self-evident.
In order to establish the second
claim, it is now enough to prove that for every
$f\in L_0(\I,X)$ the orbit map 
\[
L_0(\I,G)\ni g\mapsto g\cdot f\in L_0(\I,X)
\]
is continuous. 
Let $\e>0$ be any. 
Using the boundedness of the original 
action, choose a $W\ni e_G$ such that for all $x\in X$ and
$w\in W$, $\rho(w\cdot x, x)<\e$.
 The set
$g[W,\lambda\e/2]$ is a neighbourhood of $g$ in $L_0(\I,G)$, and
if $g_1\in g[W,\lambda \e/2]$ is arbitrary, then
for every $x\in X$ apart from a set of measure $\leq \lambda\e$
one has $\rho(g(x)\cdot f(x), g_1(x)\cdot f(x))=
\rho(f(x), w(x)\cdot f(x))<\e$, where $w(x)\equiv g(x)^{-1}g_1(x)\in W$
for every $x\in X$ apart from a set of measure $\lambda\e/2$.
This means that
$\mathrm{me}_\lambda(g_1\cdot f, g\cdot f)\leq\e$, 
establishing the continuity of the orbit map.
\end{proof}

\subsection{Proof of Theorem \ref{gener}}
Fix a parametrization of the non-atomic Lebesgue measure
space $X$, that is,
a measure space isomorphism $\I\leftrightarrow X$.
The required net of measures on $L(\I,G)$ will be indexed by the
set of all pairs of the form $(n,F)$, where $n\in\N_+$
and $F\subseteq G$ is a finite subset, directed as follows:
$(n_1,F_1)\prec (n_2,F_2)$ iff $n_1\leq n_2$ and $F_1\subseteq F_2$.
Fix a left-invariant Haar measure $\nu$ on $G$.
For every $n,F$ as above use the F\o lner condition and the assumed
amenability of the locally compact group $G$ to
choose a compact subset $K=K_{n,F}\subseteq G$ with the property
\[
\frac{\nu(gK\Delta K)}{\nu(K)}<\frac 1n
\]
for each $g\in F$.
Now let $K^n$ denote the set of all functions
in $L_0(\I,G)$ taking values in $K$ and constant on every interval
of the form $[i/n, (i+1)/n)$, $i=0,1,\dots,n-1$.
Topologically, $K^n$ can be identified
with the $n$-th power of the compact set $K$.
Denote by $\nu_{n,F}$ the
product measure $(\nu\vert_K)^n$ normalized to one and 
viewed as a probability measure on $L_0(\I,G)$ with support
$K^n$.
It remains to verify that the net of probability measures
$(\nu_{n,F})$ on the topological group $L_0(\I,G)$ satisfies
the two assumptions of Corollary \ref{levy}.
\par
(i) The net of measures $(\nu_{n,F})$ concentrates in $L_0(\I,G)$.
\par
The following general and powerful result, due to
Talagrand (\cite{Ta}, p. 76 and Prop. 2.1.1), extends the 
particular case of finite spaces belonging to Schechtman \cite{Sch,MS}.
Let $Y=(Y,\Sigma,\mu)$ denote a probability space. Then the product measures
$\mu^{\otimes n}$, $n\in\N$, on $Y^n$ concentrate, as $n\to\infty$,
with respect to the
[uniform structure generated by the] normalized 
Hamming distance on $Y^n$, given by
\[
\rho(f,g)=\frac 1n \sharp\{i\mid f_i\neq g_i\}.
\]
Moreover, the (Gaussian) bounds for the rate of concentration 
are independent of a particular $Y$, cf. {\it loco citato.}
In other words, there is a family of functions 
$\alpha_n\colon [0,1]\to[0,\frac 12]$ (of the form
$\alpha_n(\e)=C_1\exp(-C_2\e n^2)$), independent of $Y$ and $\mu$ and
such that, whenever a measurable $A\subseteq Y^n$ has the property
$\mu^{\otimes n}(A)\geq\frac 12$, one has for every $\e>0$
\[
\mu^{\otimes n}(A_\e)\geq 1-\alpha_n(\e),
\]
where $A_\e=\{y\in Y^n\colon \rho(y,A)<\e\}$.

In view of Lemma \ref{conc-trans}, it is therefore enough to 
show that the uniform structure induced on $K^n$
by the Hamming-type distance $\rho$ is finer than the
restriction of the right uniform structure $\Un_\Rsh(L(X,G))$
(which of course coincides with the unique compatible uniformity
on $K^n$).
Let $V$ be a neighbourhood of unity in $G$ and let
$\e>0$.
Let $f,g\in K^n$ be arbitrary and such that $\rho(f,g)<\e$.
Then clearly
\begin{eqnarray}
\mu(\{x\in X \mid f(x)g(x)^{-1}\notin V\}) &\leq & 
\mu(\{x\in X \mid f(x)\neq g(x)\}) \nonumber \\
&=& \frac 1n\abs{\{i\mid f_i\neq g_i\}} \nonumber \\
&=& \rho(f,g) <\e,
\end{eqnarray}
that is, $(f,g)\in [V;\e]$, establishing the claim.
\par

(ii) For every $g\in L_0(\I,G)$, the nets $(\nu_{n,F})$ and
$(g\ast\nu_{n,F})$ are asymptotically proximal.
\par
Let $g\in L_0(X,G)$.
By approximating $g$ with simple functions,
one can assume without loss in generality
that the set $F=\{g_1,\cdots,g_k\}$ of values of $g$ is finite,
and that for sufficiently large $n$, the function $g$ is constant on
each $[i/n, (i+1)/n)$. 
Since for every $g_i$ one has 
$\nu(K_{n,F}\cap g_i\cdot K_{n,F})>(1-\frac 1n)\nu(K_{n,F})$, it follows that,
whenever $n\gg k$,
\[
\nu_{n,F}(K_{n,F}^n\cap g_i\cdot K_{n,F}^n)>
\left(1-\frac 1n\right)^n\to \frac 1e.
\]
To finish the proof, notice that the restrictions of the measures
$\nu_{n,F}=(\nu\vert_K)^n$ and $g\ast \nu_{n,F}$ to 
$K_{n,F}^k\cap g_i\cdot K_{n,F}^k$ coincide, and use Remark \ref{instance}.
\qed

\section{Approximation by finite groups}

\subsection{}
The aim of this Section is to show that {\it every} (topological) group
can be approximated, albeit in a very weak sense, by finite groups.
By combining the approximation result with the extreme
amenability of Hartman--Mycielski groups, we shall later deduce 
the fixed point on compacta property for the
isometry group $\Iso(\U)$ of the complete separable Urysohn metric space.

We will state the approximation result in a few equivalent forms.
Let us say that a metric space $X$ is {\bf indexed} by a set $I$ if
there is a surjection $f_X\colon I\to X$.
We will call the pair
$(X,f_X)$ an {\bf indexed metric space.} Let us say that two metric
spaces, $X$ and $Y$, indexed with the same set $I$ are 
$\e${\bf -isometric} if for every $i,j\in I$ the distances
$d_X(f_X(i),f_X(j))$ and $d_Y(f_Y(i),f_Y(j))$ differ by at most $\e$.

\begin{lemma} If metric spaces $X$ and $Y$ indexed by a set $I$
are $\e$-isometric, then $X$ and $Y$ can be isometrically
embedded into a metric space $Z$ in such a way that for each
$i\in I$, $d_Z(f_X(i),f_Y(i))\leq\e$.
\label{gh}
\end{lemma}

\begin{proof}
Make the set-theoretic disjoint union $Z=X\cup Y$ into a weighted graph, by
joining a pair $(x,y)$ with an edge in any of the following cases:
\begin{list}{-}{}
\item $x,y\in X$, with weight $\rho_X(x,y)$;
\item $x,y\in Y$, with weight $\rho_Y(x,y)$,
\item for some $i\in I$, $x=f_X(i)$ and $y=f_Y(i)$,
with weight $\e$.
\end{list}
The weighted graph $Z$ equipped 
with the path metric clearly contains $X$ and $Y$
as metric subspaces and satisfies the required property.
\end{proof}

\begin{thm}
Let $g_1,\dots,g_n$ be a finite family of isometries of a
metric space $X$.
Then for every $\e>0$ and every finite collection
$x_1,\dots, x_m$ of elements of $X$ there exist a finite metric space
$\widetilde X$, elements $\tilde x_1,\dots,\tilde x_m$ of $\tilde X$,
and isometries $\tilde g_1,\dots,\tilde g_n$ of $\tilde X$ such that
the indexed metric spaces 
$\{g_i\cdot x_j\colon  i=1,2,\dots,n, j=1,2,\dots,m\}$ 
and $\{\tilde g_i\cdot \tilde x_j\colon  i=1,2,\dots,n, j=1,2,\dots,m\}$
are $\e$-isometric.
\label{approx}
\end{thm}

\begin{proof} 
We will perform the proof in several simple steps.
\par
1. By first rescaling the metric $\rho_X$ on $X$
and then replacing it with $\min\{\rho_X,1\}$
if necessary, we can assume without
loss in generality that the values of $\rho_X$ are bounded by $1$.
\par
2. Without loss in generality, we may also assume that $X$ supports
the structure of an abelian
group equipped with a bi-invariant metric, and $g_i$'s are
metric-preserving group automorphisms. 
For instance, one can replace $X$ with
the free abelian group $A(X)$ on $X$, and extend the metric
from $X$ to a maximal invariant metric on $A(X)$ bounded
by $1$ (the so-called {\bf Graev metric,} cf. \cite{Graev, U3}); then
every isometry of $X$ uniquely extends to an isometric
automorphism of the metric group $A(X)$.
\par
3. Let $G$ denote a group of isometries of $X$ generated
by $g_1,\dots,g_n$.
The semidirect product group $G\ltimes A(X)$ is equipped
with the bi-invariant metric $\rho$ defined by
\[
\rho((g,a),(g',a'))=\begin{cases} 
1, & \mbox{if $g\neq g'$,} \\
d(a,a'), & \mbox{otherwise.}
\end{cases}
\]
[The bi-invariance of $\rho$ is established through a 
direct calculation using the multiplication rule in
the semidirect product in question:
\[
(g,a)(h,b)=(gh,a+g\cdot b).]
\]
As usual, we will identify $G$ with a subgroup of the 
semidirect product under the
mapping $G\ni g\mapsto (g,0)$, and similarly $A(X)$ is identified with a
normal subgroup of the semidirect product under the mapping
$A(X)\ni x\mapsto (e_G,x)$.
Under such conventions,
the automorphism of $A(X)$ determined by each $g\in G$
is just $g$ itself considered as an isometric isomorphism of $A(X)$:
\begin{eqnarray*}
\forall a\in A(X), ~~gag^{-1} &\equiv & (g,0)(e,a)(g,0)^{-1} \\
&=& (g,g\cdot a)(g^{-1},0) \\
&=& (e, g\cdot a) \\
&\equiv & g\cdot a. 
\end{eqnarray*}
In particular, for every $i,j$ one has
\[
g_ix_jg_i^{-1}=g_i\cdot x_j.
\]

4. Let $F_{m+n}$ denote the free group on $m+n$ generators 
denoted by the symbols
\[
g_1,g_2,\cdots,g_n,x_1, x_2,\dots,x_m.
\]
Denote by $\pi\colon F_{m+n}\to G\ltimes A(X)$ the homomorphism
sending each generator $g_i$ to the corresponding element of $G$ and each
generator $x_j$ to the corresponding element of $X\subset A(X)$.
Pull the metric $\rho$ back from $G\ltimes A(X)$ 
to $F_{m+n}$ by letting
\[
\rho'(x,y)=\rho(\pi(x),\pi(y)).
\]
The pseudometric $\rho'$ is bi-invariant on $F_{m+n}$, though need
not be a metric.
By force of the remark at the end of step 3, 
the indexed pseudometric spaces 
\[
\{g_i\cdot x_j \mid i=1,2,\dots,n, j=1,2,\dots,m\}\subseteq
(X,\rho_X)\subset (A(X),\rho)
\] 
and
\[
\{g_i x_j g_i^{-1} \mid i=1,2,\dots,n, j=1,2,\dots,m\}\subset (F_{m+n},\rho')
\]
are isometric (and thus are both metric spaces).
\par
5. By adding to $\rho'$ an arbitrary bi-invariant metric on $F_{m+n}$
normalized so as to only slightly change the values of 
distances between pairs
of elements $g_i x_j g_i^{-1}$
(for instance, let us agree on the discrete metric
taking its values in $\{0,\e/2\})$),
we can assume without loss in generality
that $\rho'$ is a bi-invariant {\it metric} on $F_{m+n}$, while the
indexed metric spaces
$\{g_i\cdot x_j \mid i=1,2,\dots,n, j=1,2,\dots,m\}$ and
$\{g_i x_j g_i^{-1} \mid i=1,2,\dots,n, j=1,2,\dots,m\}$ are $\e/2$-isometric.

\par
6. Now replace the metric $\rho'$ with the {\bf maximal} 
among all bi-invariant metrics on $F_{m+n}$ that
coincide with $\rho'$ on the set $F^{(3)}_{m+n}$ 
of all words of reduced length $\leq 3$. 

To prove the existence of such
a metric, say $d$, denote by $\mathcal M$ the family of all bi-invariant
metrics on $F_{m+n}$ whose restriction to $F^{(3)}_{m+n}$
coincides with $\rho'$.
Since ${\mathcal M}\ni \rho'$, the family $\mathcal M$
is non-empty.
For any two elements $w,v\in F_{m+n}$ and an arbitrary
$\varsigma\in {\mathcal M}$, the value $\varsigma(w,v)$ is bounded from above
uniformly in $\varsigma$ by any sum of the form $\sum_k\rho'(a_k,b_k)$,
where $w=\sum_ka_k$ and $v=\sum_kb_k$ are two representations 
having the same length and such that $a_k,b_k\in F^{(3)}_{m+n}$. 
Now set 
\[
d(v,w)=\sup_{\varsigma\in {\mathcal M}}\varsigma (v,w).
\]
The supremum on the r.h.s. is finite, and has 
all the required properties.
\par
7. Notice that for $w,v\in F_{m+n}$ 
\[
d(w,v)=\inf\sum_k\rho'(a_k,b_k),
\]
where the infimum is taken over all possible representations of the above sort
$w=\sum_ka_k$, $v=\sum_kb_k$, having the same length and such
that $a_k,b_k\in F^{(3)}_{m+n}$. 
\par
[{\it Proof:} the infimum on the r.h.s. is a bi-invariant
pseudometric, which is greater than or equal to $d$, 
and whose restriction
to $F^{(3)}_{m+n}$ coincides with the restriction of $\rho'$.
We conclude: this infimum is in fact a metric, and it must coincide with $d$.]
\par
8. Denote by $\delta$ the smallest positive value of the distance
$\rho'$ (or, equivalently, $d$) 
between any two elements of the finite set $F^{(3)}_{m+n}$.
It follows from (7) that for every word $x\in F_{m+n}$, the value of
$d(w,e)$ is at least $[l(w)/3]\delta$, where $l(w)$ denotes the
reduced length of $w$.
\par
9. Being free, the group $F_{m+n}$ 
is residually finite, that is, admits a separating
family of homomorphisms into finite groups.
(Cf. e.g. \cite{Hall},
Ch. 7, exercise 5.)
Therefore, 
for every natural $k$ there exists a normal subgroup $N_k$
such that the factor-group $F_{m+n}/N_k$ is finite and the only
word of reduced length $\leq k$ contained in $N_k$ is identity. 
Denote by $\varpi\colon F_{m+n}\to F_{m+n}/N_k$ the factor-homomorphism
and equip $F_{m+n}/N_k$ with the factor-distance of $d$ by letting
\[
d'(g,h) = \inf\{d(w,v) \mid \varpi(w)=g, \varpi(v)=h\}.
\]
\par
10. It is immediate that $d'$ is a bi-invariant pseudometric.
Moreover, for $k\geq 3$ it is a metric:
$d'(g,h)\geq \delta\geq \e/2$ whenever $g\neq h$, cf. step 8.
\par
11. Now assume that $k\geq 3[\delta^{-1}]+4$.
Let $w,v\in F^{(3)}_{m+n}$ and $x,y\in N_k$ be arbitrary,
$x\neq y$.
Then
$d(wx,vy)=d(e,x^{-1}w^{-1}vy)\geq 1$, because $w^{-1}v\in N_k$
and therefore $l(x^{-1}w^{-1}vy)\geq 3[\delta^{-1}]$ and (8) applies.
Therefore, $d'(\varpi(w),\varpi(v))=d(w,v)$ and the restriction 
of $\varpi$ to $F^{(3)}_{m+n}$
is an isometry. 
\par
12. Now set $\tilde X=F_{m+n}/N_k$ and $\tilde x_j=\varpi(x_j)$.
For every $i=1,\dots,n$, the inner automorphism
of the finite metric group $\tilde X$ determined by $\varpi(g_i)$ is an
isometry, because the metric is bi-invariant. 
Denote this isometry by $\tilde g_j$.
The indexed metric spaces
$(g_ix_jg_i^{-1})$ and 
$(\varpi(g_i)\tilde x_j\varpi(g_i)^{-1})$,
$i=1,\dots,n, j=1,\dots,m$ are isometric by force of
the concluding remark in (11). 
Taking into account (5), we finally conclude that
the indexed metric spaces
$(\tilde g_j\cdot \tilde x_i)$ and $(g_j\cdot x_i)$ are $\e$-isometric, 
as required.
\end{proof}

\begin{remark} A careful analysis of 
the proof shows the existence of an absolute constant
$C=C(m,n,\e)>0$ such that the cardinality of the finite 
metric space in the
statement of Theorem \ref{approx} does not exceed $C$.
\endrem\end{remark}

\subsection{}
Let $G$ be a group.
One can introduce a natural topology on the
set of (equivalence classes) of all isometric actions of $G$ on metric
spaces (whose size has to be bounded from above; for instance, it
is natural to consider actions on all metric spaces $X$ of density
character not exceeding the cardinality of $G$). 
This topology is similar to the Fell topology on the set of
(equivalence classes) of unitary representations of a group
(cf. \cite{Fell} or \cite{dlHV}, p. 12), 
and is even closer to the topology
introduced by Exel and Loring on the set of representations of
a $C^\ast$-algebra \cite{EL}. 

A neighbourhood of an action $\tau$ of $G$ by isometries on a
metric space $X$ is determined by the following set of data:
a finite subset $F=\{g_1,\dots,g_n\}\subseteq G$, an $\e>0$,
and a finite collection $X'=\{x_1,\dots,x_m\}\subseteq X$.
Say that an action $\varsigma$ of $G$ on a metric space $Y$ by isometries
is in $V[F;X;\e](\tau)$, if for some finite collection
$Y'=\{y_1,\dots,y_n\}\subseteq Y$ the metric spaces
$\{g_j\cdot x_i\}$ and $\{g_j\cdot y_i\}$, naturally indexed by
$\{1,\dots,m\}\times\{1,\dots,n\}$, are $\e$-isometric. 
Our Theorem \ref{approx} can be now reformulated as follows.

\begin{corol}
Every action of a free group $F$ on an arbitrary metric space by
isometries is the limit of a net of actions of $F$ by isometries
on finite metric spaces.
\qed
\end{corol}

\begin{remark}
It is worth noting that approximation results of the above
type are not unknown.
For instance, as a corollary of a criterion by
Exel and Loring \cite{EL} and the known residual finite-dimensionality
of the group $C^\ast$-algebras of free groups \cite{GM}, every 
representation of such a $C^\ast$-algebra 
in a Hilbert space is approximated
in the Exel--Loring topology by finite-dimensional representations.
\endrem\end{remark}

To cast the above result as one on approximation of 
{\it topological} groups, we need to remind the concept of the 
Urysohn universal metric space. 

\section{Urysohn metric spaces and their groups of isometries}
We begin this Section with a summary of some known concepts and results
from theory of Urysohn metric spaces (Subsections \ref{basic}
and \ref{isom}), after which we
state a result on approximation of Polish topological groups
by finite groups (Subsection \ref{approximation}), establish
the fixed point on compacta property of the group of isometries
of the complete separable Urysohn space (Subsection \ref{thefixed}),
and finally give a new proof of the fixed point on compacta
propety for the infinite orthogonal groups (Subsection \ref{newp}).

\subsection{\label{basic} Urysohn metric spaces} 
A metric space $X$ is called a ({\bf generalized}) 
{\bf Urysohn space} if it has
the following property: whenever $A\subseteq X$ is a finite
metric subspace of $X$ and $A'=A\cup\{a\}$ is an arbitrtary one-point
metric space extension of $A$, the embedding
$A\hookrightarrow X$ extends to an isometric embedding
$A'\hookrightarrow X$.  
(Cf. \cite{Ur1,Ur,Kat,Ver,U2,GK} and \cite{Gr}, 3.11$_+$.)

There is only one, up to an isometry, complete separable Urysohn
metric space, which we denote by $\U$. 
This space contains an isometric
copy of every separable metric space.
Moreover, if $X$ is a separable
metric space and $A\subseteq X$ is a finite subspace, then every
isometric embedding $A\hookrightarrow \U$ extends to an isometric
embedding $X\hookrightarrow \U$.

A metric space $X$ is called $\mathbf n${\bf -homogeneous}, where $n$ is 
a natural number, if every isometry between two subspaces of $X$
containing at most $n$ elements each extends to an isometry of
$X$ onto itself.
If $X$ is $n$-homogeneous for every natural $n$,
then it is said to be $\omega${\bf -homogeneous}.
The complete separable Urysohn space $\U$ is $\omega$-homogeneous
and moreover enjoys the stronger property: 
every isometry between two compact subspaces of $X$ extends to
an isometry of $X$ onto itself. 
Other well-known metric spaces having the same higher
homogeneity property are the unit sphere $\s$ of the 
infinite-dimensional Hilbert space $\H$ and 
the infinite-dimensional Hilbert space $\H$ itself
(\cite{Blum}, Ch. IV, \S 38).

There are some obvious modifications of the concept of Urysohn
metric space.
For example, one can consider only metric spaces of
diameter not exceeding a given positive number $d$.
The corresponding
complete separable Urysohn space will be denoted $\U_d$.
 Another
possibility is to consider Urysohn metric spaces in
the class of metric spaces whose metrics only take values in
the lattice $\e\Z$, $\e>0$.
The corresponding object will be denoted
$\U^{\e\Z}$ (respectively, $\U_d^{\e\Z}$). \par 
Certainly, the above are not the only 
classes of metric spaces
for which the Urysohn-type universal objects exist.
For
instance, the Urysohn metric spaces for the class of spherical
metric spaces of a fixed diameter
in the sense of Blumenthal \cite{Blum}
are spheres in spaces $l_2(\Gamma)$.
The infinite-dimensional
Hilbert spaces play the role of Urysohn metric
spaces for the class of metric spaces embeddable into Hilbert spaces.

The following construction of the Urysohn space belongs to Kat\u etov 
\cite{Kat}.
Let us say, following \cite{Kat,U2,GK}, 
that a 1-Lipschitz real-valued function $f$ on a 
metric space $X$ is {\bf supported on,}
or else {\bf controlled by,} a metric subspace $Y\subseteq X$
if for every $x\in X$
\[
f(x)= \inf\{\rho(x,y) + f(y) \colon y\in Y\}.
\]
Put otherwise, $f$ is the largest 1-Lipschitz function on $X$
having the prescribed restriction to $Y$.
For instance, every distance
function $x\mapsto\rho(x,x_0)$ 
from a point $x_0$ is controlled by a singleton, $\{x_0\}$.

Let $X$ be an arbitrary metric space.
Denote by $\F(X)$ the
collection of all functions $f\colon X\to\R$ controlled by
some finite subset of $X$ (depending on the function) and
having the property
\begin{equation}
\vert f(x)-f(y)\vert\leq d_X(x,y)\leq f(x) + f(y)
\label{prop}
\end{equation}
for all $x,y\in X$.
If equipped with the supremum metric,
$\F(X)$ becomes a metric space
of the same density character as $X$,
containing an isometric copy of
$X$ under the Kuratowski embedding:
\[
X\ni x\mapsto [d_x\colon X\ni y\mapsto \rho(x,y)\in\R]\in \F(X).
\]
Besides, the space $\F(X)$ contains all one-point metric extensions of
every finite metric subspace of $X$. \par
One can form an increasing sequence of iterated extensions of the form
\[
X,\F(X), \F^2=\F(\F(X)),\dots, \F^n(X)=\F(\F^{n-1}(X)), \dots,
\] 
take the union,
$\F^\infty(X)$, and form the metric completion of it,
$\hat\F^\infty(X)$.
The latter space is a generalized Urysohn space.
If the metric space $X$ is separable, then
so is $\hat\F^\infty(X)$, and thus it is isometric to $\U$. 
If $X$ is non-separable, 
then the resulting metric space $\hat\F^\infty(X)$ 
need not be $\omega$-homogeneous.

If $X$ is a separable metric space with $\diam(X)\leq d$ and throughout 
the above construction one replaces $\F(X)$ with
the metric space $\F_d(X)$ formed by all
functions $f$ satisfying (\ref{prop}),
bounded by $d$, and controlled by finite subspaces in a suitably
modified sense, then the resulting metric space $\hat\F^\infty(X)$ 
is isometric to $\U_d$.

\subsection{\label{isom}Groups of isometries}
A remarkable feature of the above
construction, discovered by Uspenskij, is that it enables one
to keep track of groups of isometries. 
\par
Given an arbitrary metric space $X$,
the topology of pointwise convergence and the compact-open
topology on the group $\Iso(X)$ of all isometries of $X$ onto itself
coincide and turn $\Iso(X)$ into a Hausdorff topological group. 
The basic neighbourhoods of identity in this topology are of the form
\[
V[F;\e]=\{g\in\Iso(X) \colon \forall x\in F, ~ d_X(g(x),x)<\e\},
\]
where $F\subseteq X$ is finite and $\e>0$. 
If $X$ is separable (and thus second-countable), then so is $\Iso(X)$. 
\par
Notice that in general the action of $\Iso(X)$ on the metric 
space $X$ is not bounded (cf. Remark \ref{uca}.2), while
the action of $\Iso(X)$ by translations
on the space of bounded uniformly continuous (or Lipschitz)
functions on $X$, equipped with the supremum norm, is not, in
general, continuous. \par
However, the isometric action of the group
$\Iso(X)$ on the metric space of all 1-Lipschitz functions on $X$ 
controlled by finite subsets happens to be continuous.
Indeed, if
a function $f\in E(X)$ is controlled by a finite 
$Y\subseteq X$, then the translation $g\circ f$ does not differ
from $f$ by more than $\e$ at any point of $X$, provided
$g\in V[Y;\e]$.
Consequently, the canonical 
representation of $\Iso(X)$ in $\F(X)$ by isometries
defines a topological group embedding
$\Iso(X)\hookrightarrow\Iso(\F(X))$. \par
Iterating this process countably many times, one obtains a
a continuous action of $\Iso(X)$ by
isometries on $\F^\infty(X)$, which in its turn
extends to a continuous action of $\Iso(X)$ on
the metric completion $\hat\F^\infty(X)\cong\U$. \par
We adopt terminology suggested in \cite{U2}
and say that a metric subspace $Y$ 
is {\bf $g$-embedded} into a metric space $X$ if there exists 
an embedding of topological groups
$e\colon\Iso(Y)\hookrightarrow\Iso(X)$ with the property that
for every $h\in\Iso(Y)$ the isometry $e(h)\colon X\to X$ is
an extension of $h$. 
The above argument establishes the following result.

\begin{prop}[Uspenskij \cite{U1}]
Every separable metric space $X$ can be $g$-embedded into 
the complete separable Urysohn metric space $\U$.
\label{g-emb}
\qed
\end{prop}

Since every [second-countable] topological group $G$ embeds into
the isometry group of a suitable [separable] metric space
\cite{Te}, we arrive at the following.

\begin{thm}[Uspenskij \cite{U1}]
The topological group $\Iso(\U)$ is the universal second-countable
topological group.  \qed
\label{univ}
\end{thm}

(Cf. also \cite{Gr}, 3.11.$\frac 23_+$.)

Since every isometry between two compact subspaces of $\U$
can be extended to an isometry of $\U$ onto itself, we obtain
the following useful corollary of Proposition \ref{g-emb}.

\begin{corol}
\label{conven}
Each isometric embedding of a compact metric space into $\U$ is a
$g$-embedding.
\qed
\end{corol}

The question of the existence of universal topological groups of a
given uncountable weight $\tau$  
(in fact, of {\it any} uncountable weight $\tau$) 
remains open. 
However, recently Uspenskij has established the following result. 

\begin{thm}[\cite{U2}]
Every topological group $G$ embeds, as a topological subgroup,
into the group of isometries $\Iso(X)$ of a suitable $\omega$-homogeneous
Urysohn metric space $X$ of the same weight as $G$.
\label{uspenskij}
\end{thm}

The construction rather resembles the proof of Theorem
\ref{univ}, but in order to achieve $\omega$-homogeneity of the
union space, one alternates between
the Kat\u etov metric extension $\F(\cdot)$
and the `homogenization' extension, $H(\cdot)$,
which forms the nontrivial technical core of the proof and is
described in the following theorem.

\begin{thm}[Uspenskij \cite{U2}]
Every metric space $X$ $g$-embeds into an $\omega$-homogeneous metric
space $H(X)$ of the same weight as $X$. \qed
\label{homo}
\end{thm}

\subsection{\label{approximation}
Approximation of topological groups}

Now we can state yet another 
reformulation of the approximation Theorem \ref{approx}.

\begin{thm} 
For every finite collection of isometries $g_1,\cdots,g_n$ of
the complete separable Urysohn metric space
$\U$ and every neighbourhood $V$ of identity in $\Iso(\U)$
there are isometries $h_1,\dots,h_n\in\Iso(\U)$ generating a
finite subgroup and such that $h_ig_i^{-1}\in V$, $i=1,\cdots,n$.
\label{approxurysohn}
\end{thm}

\begin{proof} One can assume that $V=V[X;\e]$,
where $X=\{x_1,\cdots,x_m\}\subseteq \U$ and $\e>0$. 
Using Theorem \ref{approx}, choose a finite metric space
$\widetilde X$, elements $\tilde x_1,\dots,\tilde x_m$ of $\tilde X$,
and isometries $\tilde g_1,\dots,\tilde g_n$ of $\tilde X$ such that
the naturally indexed finite metric spaces 
$A=\{g_i^{-1}\cdot x_j \mid i=1,2,\dots,n, j=1,2,\dots,m\}$ and
$B=\{\tilde g_i^{-1}\cdot \tilde x_j \mid i=1,2,\dots,n, j=1,2,\dots,m\}$
are $\e/2$-isometric.
\par Using Lemma \ref{gh}, isometrically
embed $A$ and $\tilde X$ into a finite metric space
$Z$ in such a way that $d_Z(g_i^{-1}(x_i),\tilde g_i^{-1}(\tilde x_i))
\leq \e/2$ for all $i,j$.
Now extend the embedding $A\hookrightarrow \U$ to an isometric
embedding $Z\hookrightarrow \U$.
According to Corollary \ref{conven},
the (finite) group $\Iso(\tilde X)$ simultaneously extends to a
group of isometries of $\U$.
Denote the extension of the
isometry $\tilde g_i$ by $h_i$.
One has for all $i,j$:
\[
d(x_i,h_jg_j^{-1}(x_i))=d(\tilde g_i^{-1}(x_i),g_i^{-1}(x_i))
\leq\e/2<\e,
\] 
and the proof is finished.
\end{proof}

Let $G$ be a group and let $X$ be a metric space. 
Every action of $G$ on $X$ by isometries can be viewed as a homomorphism 
$\tau\colon G\to\Iso(X)$.
Equip the set $\mathrm{Hom}(G,\Iso(X))$ of
all such homomorphisms with the topology of pointwise convergence
on $G$, that is, the one induced from the Tychonoff product 
$\Iso(X)^G$. Since, in its turn, the topological space $\Iso(X)$
is a subspace of the Tychonoff product $X^X$, one concludes that
$\mathrm{Hom}(G,\Iso(X))$ is a topological subspace of the Tychonoff
product $X^{G\times X}$. In this form, the identification of 
the collection of all actions
$\tau\colon G\times X\to X$ with a subspace of $X^{G\times X}$ 
becomes obvious.

Call an action {\bf periodic} if it factors through an
action of a finite group. One can reformulate 
Theorem \ref{approxurysohn} as follows.

\begin{corol} Let $F$ be a free group.
The set of 
periodic actions of $F$ on the Urysohn metric space $\U$ is
everywhere dense in the set of all actions. \qed
\end{corol}

Let $F_\infty$ denote the free group of countably infinite rank.
The mapping associating to an action $\tau$ of $F_\infty$ on $\U$
the closure of $\tau(F_\infty)$ in $\Iso(\U)$ is a surjection
from $\mathrm{Hom}(G,\Iso(X))$ onto the space ${\mathcal L}(\Iso(\U))$
of all closed subgroups of $\Iso(\U)$.
Equip the latter space with the corresponding quotient topology.
The topology so defined satisfies the axiom $T_0$. 

\begin{corol}
The set of finite subgroups is everywhere dense in the 
topological space ${\mathcal L}(\Iso(\U))$. \qed
\end{corol}

This leads to an approximation result for Polish topological
groups.

\begin{corol}
Let $G$ be a Polish topological group. 
Then under every isomorphic
embedding into $\Iso(\U)$ the group $G$ is the limit of a net of
finite subgroups. \qed
\label{polish}
\end{corol}

\begin{remark}
At the first sight, the above may seem to contradict the
general principle (in particular espoused and explained by Vershik 
in \cite{Ver1}) according to which approximability of an (infinite)
group $G$ by 
finite groups is essentially equivalent to amenability of $G$. 
In fact, our results are in
perfect agreement with this principle in that the approximating
groups come from `without' the group $G$ and thus form
an approximation not to $G$ itself, but to a suitable
topological group extension of $G$, which indeed turns out to be 
amenable (and even extremely amenable).
\endrem\end{remark}

\subsection{\label{thefixed}
The fixed point property of the group $\Iso(\U)$}

Theorems \ref{approxurysohn} and \ref{gener} enable us
to deduce the fixed point on compacta property for
the group of isometries of the complete separable Urysohn space $\U$.

\begin{thm} The group $\Iso(\U)$ of all isometries of the
complete separable Urysohn space $\U$, equipped with the
standard {\rm (}pointwise $=$ compact-open{\rm )} topology, is extremely
amenable {\rm (}has the fixed point on compacta property{\rm )}.
\label{fpc}
\end{thm}

\begin{proof}
Let the group $\Iso(\U)$ act continuously on a compact space $K$. 
We will show that every finite collection of elements of $\Iso(\U)$
has a common fixed point in $K$, from which the result follows by an
obvious compactness argument. 
Fix an arbitrary such collection, $g_1,\dots,g_n\in\Iso(\U)$.

Let $U\in\Un_K$ be an arbitrary element of
the unique compatible uniform structure on $K$.
Without loss in generality, assume that $U$ is closed as a 
subset of $K\times K$ (and consequently compact).
Using the boundedness of the action of $\Iso(X)$ on $K$, choose a finite
$X\subseteq\U$ and an $\e>0$ such that, whenever $g\in V=V[X;\e]$,
one has $(g\cdot\kappa, \kappa)\in V$ for all $\kappa\in K$.

By Theorem \ref{approxurysohn}, 
there are isometries $h_1,\dots,h_n\in\Iso(\U)$ generating a
finite subgroup $H$ and such that $h_ig_i^{-1}\in V$, $i=1,\cdots,n$.

Let $\tilde X$ be a finite $H$-invariant subset of $\U$ containing $X$.
The iterated Kat\u etov extension
$\hat\F^\infty(L_0(\I,\tilde X))$ contains $\tilde X$ as a subspace made up
of all constant functions and is isometric to $\U$, and
since $\tilde X$ is finite, an isometry between the two spaces can
be chosen so as to extend the canonical embedding of $\tilde X$ into $\U$.
Thus we obtain a chain of $g$-embeddings
\[
\tilde X\subset L_0(\I,\tilde X) \subset \U.
\]
The group $L_0(\I,H)$ acts on $L_0(\I,\tilde X)$  
continuously and isometrically (Lemma \ref{L}), and this action canonically
extends to a continuous isometric action of the same group on
the space $\hat\F^\infty(L_0(\I,\tilde X))\cong\U$.
Thus we obtain a
continuous group monomorphism $j\colon L_0(\I,H)\to\Iso(\U)$ with
the property that
for every $h\in H$ one has $j(h)\vert_{\tilde X} = h\vert_{\tilde X}$.  \par
Composing $j$ with the action $\Iso(\U)\to\Homeo(K)$, we obtain a
continuous action of $L_0(\I,H)$ on $K$. 
By force of Theorem \ref{gener}, $L_0(\I,H)$ has a 
common fixed point in $K$, say $\kappa$.
In particular, $\kappa$ is fixed under the elements 
$j(h_1),\dots,j(h_n)\in\Iso(\U)$, where we identify elements of $H$
with constants in $L_0(\I,H)$. \par
For all $x\in X$ and $i=1,2,\dots,n$, one has
$d(j(h_i)^{-1}(x),g_i^{-1}(x))=d(h_i^{-1}(x),g_i^{-1}(x))<\e$
for all $i$ and $x\in X$, implying that 
$j(h_i)g_i^{-1}\in V$ for $i=1,2,\dots,n$. 
Consequently and by the choice of $V=V[X;\e]$,
\[
(g_i\kappa,\kappa)\equiv(g_i\kappa,j(h_i)\kappa)\equiv 
(g_i\kappa,(j(h_i)g_i^{-1})\cdot(g_i\kappa))\in U
\] 
for all $i$.
Denote by $F_U$ the (non-empty) set of all
points $x\in K$ with the property $(g_ix,x)\in U$ for all $i$.
Since $U$ is closed, so is $F_U\subseteq K$.
If $U_1\subseteq U_2$,
then $F_{U_1}\subseteq F_{U_2}$.
It means that $\{F_U\}$ is
a centred system of closed subsets of the compact space $K$
and therefore has a common point, which is clearly fixed
under $g_1,\dots,g_n$, as required.
\end{proof}

\begin{remark}
The same argument {\it verbatim}
also establishes the fixed point on compacta
property of the topological group $\Iso(\U_d)$ of isometries of
the complete separable universal Urysohn space of finite diameter $d$.
\endrem\end{remark}

\subsection{\label{newp}
A new proof of the fixed point on compacta property 
of the infinite orthogonal group}
The above proof can be easily modified so as to result in a new proof
of extreme amenability of the orthogonal group $\mathrm{O}(\H)$
of an infinite-dimensional Hilbert space with the strong
operator topology.
This proof does not rely on such advanced tools from geometry as Gromov's 
isoperimetric inequality for groups $\mathrm{SO}(n)$.

The following belongs to folklore.

\begin{lemma}
Let $X$ be a metric subspace of the unit sphere $\s$ of a real
Hilbert space $\H$.
Suppose a topological group $G$ acts on
$X$ continuously by isometries.
Then the action of $G$ extends
to a strongly continuous action of $G$ by isometries on the
sphere $\s$ {\rm (}that is, to a strongly continuous orthogonal
representation of $G$ in $\H${\rm )}. 
Put otherwise, every metric subspace of the unit sphere $\s$ of a
real Hilbert space is $g$-embedded into $\s$.
If the linear span of $X$ is dense in $\H$, the
extension is unique. 
\label{extend}
\end{lemma}

\begin{proof}
Since for every $x,y\in X$ the value of the inner product is
uniquely determined by the Euclidean distance between the elements,
\[
(x,y)=1-\frac 12 \rho_X(x,y)^2,
\]
there is only one way to turn
the linear span $\mathrm{lin}(X)$ into a pre-Hilbert
space so as to induce the given metric on $X$.
The corresponding completion 
${\mathcal K}=\widehat{\mathrm{lin}}(X)$ 
is isometrically isomorphic to the closed linear span of $X$ in
$\H$, that is, $\H=\mathcal{K}\oplus\mathcal{K}^\perp$.
As another consequence of the same observation, 
every isometry of $X$ lifts to a unique
orthogonal transformation of $\mathcal{K}$. 
The resulting homomorphism
$\pi\colon G\to\mathrm{O}(\mathcal{K})$ is 
continuous if the latter group is
equipped with the topology of simple convergence on $X$ or, which is
the same, on $\mathrm{lin}(X)$.
On the groups of isometries of metric spaces
the topology of simple convergence on an everywhere dense subset
coincides with the topology of simple convergence
on the entire space.
Consequently, the extended orthogonal representation $\pi$ of $G$ in 
$\mathcal K$ is strongly continuous.
It remains to extend $\pi$
to a representation $\left(\begin{matrix}\pi & 0 \\
0 & \mathrm{Id}_{\mathcal{K}^\perp}\end{matrix}\right)$ of $G$ in $\H$.
The uniqueness statement is obvious.
\end{proof}

Here is an outline of the alternative 
proof of extreme amenability of
$\mathrm{O}(\H)_s$.
We will be only considering the separable case
$\H=l_2$; the extension to non-separable case is straightforward.
\par
Every finite collection $g_1,g_2,\dots,g_n$ of
elements of $\mathrm{O}(l_2)$, viewed as isometries of the unit sphere
$\s$, can be approximated (in the strong operator topology) 
by a collection of elements $g'_1,g'_2,\dots,g'_n$ of
a finite-dimensional orthogonal subgroup in the following 
sense: for a given natural $m$ and an $\e>0$, one has
$\norm{g_i(e_j)-g'_i(e_j)}<\e$ for all $i=1,2,\dots,n$, 
$j=1,2,\dots,m$, where $g'_j\in \mathrm{O}(N)$,
$e_j$ denote the standard basic vectors, and the rank
$N$ is sufficiently large.   \par
According to Lemma \ref{L}, the topological group $L_0(\I,\mathrm{O}(N))$  
acts continuously by isometries
on the metric space $L_2(\I,\s^N)$, equipped with the $l_2$-metric.
(The topology induced on $L(\I,\s^N)$ by $l_2$-metric
is still that of convergence in measure, because $\s^N$ is compact.)
The metric space $L_2(\I,\s^N)$ 
is spherical of diameter one and thus can be embedded into $\s$ as a metric
superspace of $\s^N$.
Using Lemma \ref{extend}, we obtain a chain of continuous monomorphisms of
topological groups
\[
\mathrm{O}(N)<L_0(\I,\mathrm{O}(N))<\Iso(L_2(\I,\s^N))
< \mathrm{O}(l_2).
\]
According to Theorem \ref{gener},
the second topological group on the left is extremely amenable.
It follows
that the orthogonal operators $g'_1,g'_2,\dots,g'_n$ have a common
fixed point in every compact space upon which $\mathrm{O}(l_2)$
acts continuously.
Now the proof is accomplished in the same way
as in Theorem \ref{fpc}. 

\section{Ramsey-type theorems for metric spaces vs f.p.c. property}

\subsection{Ramsey--Dvoretzky--Milman property}
In order to extend the result about fixed point on compacta
property of the isometry group $\Iso(\U)$ beyond 
the separable case, we will obtain a new characterization of
extremely amenable groups of isometries in terms of a Ramsey-type
property of the metric spaces $X$.

The following is an adaptation from \cite{Gr1}, Sect. 9.3.

\begin{defin}
Let $G$ be a group of uniform isomorphisms of a uniform 
space $X$.
We will say that the pair $(G,X)$ has the 
{\bf Ramsey--Dvoretzky--Milman property} if 
for every bounded uniformly continuous function
$f$ from $X$ to a finite-dimensional Euclidean space, 
every $\e>0$, and every
compact $K\subseteq X$, the function $f$ is $\e$-constant on a suitable
translate of $K$, that is, there is a $g\in G$ such that 
\[
\mathrm{Osc}(f\mid gK)<\e.
\]
Equivalently, `compact' can be replaced with `finite.'
\endrem\end{defin}

We defer two master examples (Ex. \ref{secmaster}
and \ref{milmaster}) in order to precede them by a few simple
preliminary results.
The following is established by pulling back the function
$f$ from $Y$ to $X$.

\begin{lemma}
Let $G$ be a group, acting by uniform isomorphisms on the uniform
spaces $X$ and $Y$, and let $f\colon X\to Y$ be an equivariant
uniformly continuous map with everywhere dense range.
If the pair $(G,X)$ has the Ramsey--Dvoretzky--Milman property, 
then so does $(G,Y)$.
\qed
\label{pulling}
\end{lemma}

Denote by $\Un^\ast_X$ the totally bounded replica of the
uniform structure $\Un_X$ on $X$, that is, the coarsest uniform structure
preserving the uniform continuity of every 
bounded uniformly continuous function on $X$.
Basic entourages of the diagonal for 
$\Un^\ast_X$ are of the form
\[
\{(g,h)\in X\times X \colon \abs{f(x)-f(y)}<\e\},
\]
where $f\colon X\to\R^N$ is bounded uniformly continuous, $N\in\N$.

The following reformulation of the R--D--M property is immediate.

\begin{prop}
A pair $(G,X)$ has the Ramsey--Dvoretzky--Milman property if and only if
for every compact (equivalently: finite) 
$K\subseteq X$ and every entourage
$V\in\Un^\ast_X$ there is a $g\in G$ with $gK$ being $V$-small:
$gK\times gK\subseteq V$.
\qed
\end{prop}

\begin{prop} Let $X=(X,\Un_X)$ be a uniform space.
A basis of entourages for the totally bounded replica $\Un^\ast_X$
of $\Un_X$ is given by all finite covers of the form
$\{V[A]\colon A\in\gamma\}$, where $\gamma$ is an arbitrary finite
cover of $X$ and $V\in\Un_X$.
\end{prop}

\begin{proof} The claim consists of two parts: first,
that all sets of the form
\[
\cup_{A\in\gamma} V[A]\times V[A], \mbox{ $\gamma$ finite, } V\in\Un_X
\]
are elements of $\Un^\ast_X$, and second, that each enourage from
$\Un^\ast_X$ contains a set of the above type. 

(1) Given $\gamma$, $V$, and $A$ as above, choose a bounded uniformly
continuous pseudometric $d$ on $X$ such that $(d(x,y)<1)
\Rightarrow ((x,y)\in V)$, and introduce a bounded uniformly continuous 
function $f$ from $X$ to the Euclidean space $\R^{\abs\gamma}$ with
each component $f_A$, $A\in\gamma$, defined by
\[
X\ni x\mapsto f_A(x):=d(x,A)\in\R.
\]
The set $\{(x,y)\in X^2\colon \abs{f(x)-f(y)}<1\}$ is an element of
$\Un^\ast_X$ and a subset of $\cup_{A\in\gamma} V[A]\times V[A]$.

(2) Let $W\in \Un^\ast_X$ be arbitrary. Choose a bounded
uniformly continuous function
$f\colon X\to\R^N$ and an $\e>0$ such that
$\{(x,y)\in X^2\colon \abs{f(x)-f(y)}<\e\}\subseteq W$.
Partition the image $f(X)$ into finitely many
pieces of diameter $\leq\e/2$ each and
let $\gamma$ be the family of preimages of those pieces under $f$. 
Define $V=\{(x,y)\in X^2\colon \abs{f(x)-f(y)}<\e/2\}
\in \Un^\ast_X\subseteq \Un_X$. Clearly,
$\cup_{A\in\gamma} V[A]\times V[A]\subseteq W$.
\end{proof} 

As an immediate corollary, one obtains the following.

\begin{prop}
A pair $(G,X)$ has the Ramsey--Dvoretzky--Milman property if and only if
for every compact (equivalently: finite) 
$K\subseteq X$, every finite cover $\gamma$ of $X$, and every entourage
$V\in\Un_X$, there is a $g\in G$ such that $gK$ is contained in the 
$V$-neighbourhood of some $A\in\gamma$. 
\label{cover}
\qed
\end{prop}

Here is the first major example.

\begin{example}
\label{secmaster}
Let $\Gamma$ be an infinite set, and let $n$ be a natural number. 
Choose as $G$ the group $S^f_\Gamma$ of all finite permutations of $\Gamma$,
and as $X$ the set $\Gamma^{(n)}$ of all $n$-subsets of $\Gamma$, equipped
with the finest (discrete) uniformity. 
Using Proposition \ref{cover}, one can easily see that 
the pair $(\Gamma^{(n)},S^f_\Gamma)$ has the 
Ramsey--Dvoretzky--Milman property, which statement is indeed equivalent to 
the finite Ramsey theorem.
\endrem
\end{example}

Recall that the basic entourages for the left uniform structure
$\Un_\Lsh(G)$ on a topological group $G$ are of the form
\[
V_\Lsh=\{(g,h)\in G \times G \colon g^{-1}h\in V\},
\]
where $V$ is a neighbourhood of identity in $G$. 
If $d$ is a left invariant continuous pseudometric on $G$
and $\e>0$, then the set $V[d;\e]=\{(x,y)\in X^2\colon d(x,y)<\e\}$ is
an element of $\Un_\Lsh(G)$. 
Since for every neighbourhood of identity $V$ there is a
bounded left invariant continuous pseudometric $d$ on $G$ with
$(d(x,e_G)<1)\Rightarrow (x\in V)$ and consequently
$V_\Lsh\supseteq V[d;1]$, it follows that the left uniform structure on
a topological group is determined by left invariant bounded
continuous pseudometrics.

If $d$ is a left invariant continuous pseudometric on $G$, then
$H_d=\{x\in G\colon d(x,e_G)=0\}$ forms a closed subgroup of $G$,
and the pseudometric $d$ induces a continuous left-invariant metric
$\hat d$ on the factor-space $G/H_d$ by the formula
$\hat d(xH,yH):=d(x,y)$. The canonical factor-map 
$\pi\colon G\to (G/H_d,\hat d)$ is uniformly continuous.
Notice that in general both the topology and the uniform
structure induced by $\hat d$ are coarser than
the factor-topology and the left uniform structure 
on $G/H_d$. We will denote the $G$-space $G/H_d$ equipped with the
left invariant metric $\hat d$ by $G/d$, which is consistent with the
notation sometimes used in set-theoretic topology: in our situation, $G/d$ 
is the metric space canonically associated to the pseudometric space
$(G,d)$. 

The following result (which grew out of V.V. Uspenskij's conjecture)
reveals the link between the Ramsey--Dvoretzky--Milman property and 
the existence of fixed points.

\begin{thm}
For a topological group $G$, the following are equivalent.
\begin{enumerate}
\item[(i)] $G$ has the fixed point on compacta
property.
\item[(ii)] The pair $(G,G_\Lsh)$ has the Ramsey--Dvoretzky--Milman property.
\item[(iii)] For every left-invariant continuous pseudometric $d$ on
$G$, the pair $(G,G/d)$ has the Ramsey--Dvoretzky--Milman property.
\item[(iv)] Whenever $G$ acts continuously and transitively
by isometries on a metric
space $X$, the pair $(G,X)$ has the Ramsey--Dvoretzky--Milman property.
\item[(v)] For some family $D$ of bounded continuous left invariant
pseudometrics $d$, generating the topology of $G$, each 
pair $(G, G/d)$ has the Ramsey--Dvoretzky--Milman property.
\end{enumerate}
\label{rdm}
\end{thm}

\begin{proof} (i) $\Leftrightarrow$ (ii): according to Theorem
\ref{folk}, the fixed point on compacta property of a topological
group $G$ is equivalent to the following: for every bounded
{\it right} uniformly continuous function $f$ on $G$ taking values
in a finite-dimensional Euclidean space, every finite collection
of elements $g_1,g_2,\dots,g_n\in G$, and every $\e>0$, 
there is an $x\in G$ such that
\[
\abs{f(x)-f(g_ix)}<\e \mbox{ for all } i=1,2,\dots,n.
\]
The mirror image of the above statement applies to {\it left}
uniformly continuous functions and calls for the existence of
an $x\in G$ with the property
$\abs{f(x)-f(xg_i)}<\e$ for all $i$.
This amounts to the
Ramsey--Dvoretzky--Milman property for the pair $(G,G_\Lsh)$ 
relative to the left action (with $K=\{e_G,g_1,g_2,\dots,g_n)$).
\par
(ii) $\Rightarrow$ (iii): as the canonical map $G\to G/d$ is
uniformly continuous and $G$-equivariant,  Lemma \ref{pulling} applies.
\par
(iii) $\Rightarrow$ (iv):
Let $d_X$ denote the invariant metric on $X$. 
Fix an arbitrary point $x_0\in X$. The
formula $d(g,h):=d_X(gx_0,hx_0)$ defines a left invariant continuous 
pseudometric on $G$, and the map $G\ni g\mapsto gx_0\in X$ factors through to
a $G$-equivariant isometric isomorphism between $G/d$ and $X$. 
\par
(iv) $\Rightarrow$ (v): Trivial, as $G$ acts on each space $G/d$ 
continuously and transitively by isometries.
\par
(v) $\Rightarrow$ (ii): Suppose we are given a finite subset 
$F\subseteq G$, a finite cover $\gamma$ of $G$, and a basic element 
$V_\Lsh$ of the left uniformity $G_\Lsh$, where $V$ is a
neighbourhood of identity in $G$. 
Choose a bounded left invariant continuous pseudometric $d\in D$
with the property $(d(x,e_G)<1)\Rightarrow (x\in V)$.
The sets $\pi(A)$, $A\in\gamma$, where $\pi\colon G\to G/d$ is the
factor-map, form a finite cover of $G/d$, and by assumption there
is a $g\in G$ such that $g\pi(F)$ is entirely contained in the
$1$-neighbourhood of some $\pi(A)$, $A\in\gamma$. 
Denote, as before, $H_d=\{g\in G\colon d(g,e_G)=0\}$. 
The value of $d$ is independent on the choice of representatives in
left cosets: $d(xh_1,yh_2)=d(x,y)$ for all $x,y\in G$, $h_1,h_2\in H_d$.
Let $f\in F$ be any. Since $\hat d(g\pi(f),\pi(a))<1$
for some $a\in A$, one has $d(gf,a)<1$,
that is, $gF$ is contained in $V_\Lsh[A]=AV$, 
and the Ramsey--Dvoretzky--Milman property of $(G,G_\Lsh)$ 
is thus verified.
\end{proof}

\begin{example}
\label{milmaster} The second major example is given by the pair
consisting of the full
unitary group $U(\H)$ of an infinite-dimensional Hilbert space
$\H$ and the unit sphere $\s_\H$ equipped with the Euclidean
distance. The Ramsey--Dvoretzky--Milman
property of this pair follows from Th. \ref{rdm} and the extreme
amenability of $U(\H)_s$ (cf. Subsection \ref{newp}). 
In fact, a direct proof of this property does not require
the extreme amenability of the unitary group, and such was the
original proof by Milman \cite{M3} (who then used the R--D--M property 
to give a new proof
of Dvoretzky theorem on almost spherical sections of convex bodies),
cf. also \cite{Gr1}, Sect. 9.3.
\end{example}

For sufficiently homogeneous spaces and their full groups of isometries
Theorem \ref{rdm} assumes a combinatorial form of a Ramsey-type
result for metric spaces somewhat in the spirit of 
\cite{Mat} or \cite{Kom}, but in an `approximate' implementation. 
We proceed to examine this connection now.

\subsection{Ramsey-type properties of metric spaces}
Let $X$ be a metric space, and let $F$ be a finite metric subspace of $X$.
The stabilizer of $F$,
\[
\St_F=\{g\in\Iso(X)\colon gx=x\mbox{ for each } x\in F\},
\]
is a closed subgroup of $\Iso(X)$.
Denote by  $X^{\hookleftarrow F}$ the family of all isometric embeddings
of $F$ into $X$, equipped with the natural 
action of $\Iso(X)$ on the left:
\[
X^{\hookleftarrow F}\ni j\mapsto g\circ j \in X^{\hookleftarrow F}.
\]
The supremum metric on $X^{\hookleftarrow F}$, given by
\[
d_{sup}(i,j)=\max\{d(i(x),j(x)) \colon x\in F\},
\]
is $\Iso(X)$-invariant. Denote by $d_F$ the pull-back of the metric 
$d_{sup}$ to $\Iso(X)$:
\[
d_F(g,h):=d_{sup}(gi_F,hi_F),
\]
where $i_F\colon F\hookrightarrow X$ is the canonical embedding. 
Left-invariant pseudometrics
of the form $d_F$, where $F$ runs over all finite subspaces of $X$,
generate the usual topology of pointwise convergence on $\Iso(X)$. 

If $X$ is $\abs F$-homogeneous,  
the following establishes an isomorphism of $\Iso(X)$-sets:
\begin{equation}
\label{isomo}
G/\St_F\ni g\St_F \mapsto [g\vert_F \colon F \to gF]\in 
X^{\hookleftarrow F}.
\end{equation}

In the combinatorial spirit, we 
will refer to [finite] partitions of a metric space $X$ as
{\bf colourings} of $X$ [using finitely many colours]. 
A subset $Y\subseteq X$ is {\bf monochromatic} if $Y\subseteq X$
for some $A\in\gamma$, and {\bf monochromatic up to an $\e>0$}
if $Y$ is contained in the $\e$-neighbourhood of some $A\in\gamma$.

A direct application of Theorem \ref{rdm} now results in the following.

\begin{thm} Let $X$ be an $\omega$-homogeneous 
metric space.
The following conditions are equivalent.

\begin{enumerate}
\item[(i)] The full group of
isometries $\Iso(X)$ with the pointwise topology is extremely
amenable.
\item[(ii)] Let $F\subseteq X$ be a finite metric space, and let
$X^{\hookleftarrow F}$ be coloured using finitely many colours.
Then for every finite metric subspace $G\subseteq X$ and every $\e>0$
there is an isometric copy of $G$, $G'\subseteq X$, such that
all isometric embeddings $F\hookrightarrow X$ that factor through
$G'$ are monochromatic up to $\e$.
\end{enumerate}
\label{ramsey}
\qed
\end{thm}

\begin{remark}
There is a natural surjection from
$X^{\hookleftarrow F}$ onto the collection $X^{(F)}$ of all subspaces
of $X$ isometric to $F$, as the latter space is obtained
from the former one by factoring out the group of distance-preserving
permutations of $F$:
\[
X^{(F)}\cong X^{\hookleftarrow F}/\Iso(F).
\]
In particular, if the metric space $F$ is rigid (for example, if
no two distances between different pairs of points are the same),
then the spaces $X^{(F)}$ and $X^{\hookleftarrow F}$ can be
identified.
In general, however, the distinction between the
two spaces has to be maintained, and as we shall see (Theorem
\ref{doesnot}), some groups
of isometries of $\omega$-homogeneous metric spaces
fail to have the fixed point on compacta property 
namely due to the fact that the two
spaces $X^{(F)}$ and $X^{\hookleftarrow F}$ are different.
\endrem\end{remark}

\begin{remark}
Theorem \ref{ramsey} provides at one's disposal
a rather versatile tool.
The main application in this article will
be to establish the extreme amenability of groups of isometries
of $\omega$-homogeneous generalized Urysohn spaces.
The result shall also be used to demonstrate 
that some groups of isometries are not extremely amenable.
And finally, one can turn Theorem \ref{ramsey} around in order to
deduce Ramsey-type results for metric spaces from the
known results on extreme amenability of various 
topological groups established
by other means.
The next Section contains examples of applications
of each sort.
\endrem\end{remark}

\section{Applications\label{secappl}}
\subsection{Extreme amenability of the groups $\Iso(U)$.} 
We want to formalize the content of the condition (ii)
of Theorem \ref{ramsey}, as follows.

\begin{defin}
Let $F$ and $G$ be finite metric spaces, 
let $m\in\N_+$, and let $\e>0$.
Denote by $R(F,G,m,\e)$ the following property of a metric space $X$:
\par
$X\in R(F,G,m,\e)$ $\Leftrightarrow$ {\it for every 
colouring of the set $X^{\hookleftarrow F}$ of all isometric embeddings
of $F$ into $X$ with $\leq m$ colours, there is
an isometric embedding $j\colon G\hookrightarrow X$
such that all embeddings
of $F$ into $X$ that factor through $j$ are monochromatic up to $\e$.}
\par
Say that a metric space $X$ has property $R$ if $X\in R(F,G,m,\e)$
for all finite metric spaces $F,G$ embeddable into $X$,
for all $m\in\N$, and all $\e>0$.
\endrem\end{defin}

\begin{remark}
Now Theorem \ref{ramsey} can be reformulated as
follows: {\it an $\omega$-transitive 
metric space $X$ has property $R$ if
and only if the topological group $\Iso(X)$ is extremely amenable.}
\endrem\end{remark}

\begin{prop}
Let  $F$ and $G$ be finite metric spaces, let $X$ be a metric space
containing a copy of $F$, let $m$ be a natural number, and let $\e>0$.
The following are equivalent.
\begin{enumerate}
\item[(i)] $X\in R(F,G,m,\e)$.
\item[(ii)] There is a
finite subspace $Z\subseteq X$ containing a copy of $F$ 
such that $Z\in R(F,G,m,\e)$.
\end{enumerate}
\end{prop}

\begin{proof}
(i) $\Rightarrow$ (ii): assume $\neg$(ii), that is,
no finite subspace $Z$ of $X$ containing a copy of $F$ is in $R(F,G,m,\e)$.
Denote by $\mathcal Z$ the collection of all finite metric subspaces
$Z\subseteq X$ with $Z^{\hookleftarrow F}\neq\emptyset$. By assumption,
${\mathcal Z}\neq\emptyset$.

Then for every $Z\in {\mathcal Z}$ the set $Z^{\hookleftarrow F}$ 
admits a colouring with $m$ colours, which we will view as a function
$f_Z\colon Z^{\hookleftarrow F}\to \{1,2,\dots,m\}$, 
in such a way that the following holds:
\par
($\ast$) for every isometric embedding $i\colon G\hookrightarrow Z$
and every colour $k=1,2,\dots,m$ there is an isometric embedding
$j_k\colon F\hookrightarrow G$ such that the $\e$-neighbourhood of
$i\circ j_k$ in $Z^{\hookleftarrow F}$ contains no elements of
colour $k$.
\par 
The system $\mathcal Z$ is directed by inclusion, and the collection
of intervals $[K,\infty)=\{Z\in {\mathcal Z} \colon K \subseteq Z\}$,
where $K\subseteq X$ is finite, is a filter on $\mathcal Z$,
which we will denote by $\mathcal F$.
Since $X$ can be assumed infinite (otherwise there is nothing to
prove), $\mathcal F$ 
extends to a free ultrafilter $\Lambda$ on $\mathcal Z$.
For every $j\in X^{\hookleftarrow F}$, 
one has $[\{j(F)\},\infty)\in{\mathcal F}\subset
\Lambda$, and therefore exactly one
of the sets $\{Z\in {\mathcal Z} \colon f_Z(j)=i\}$, 
$1\leq i\leq m$ is in $\Lambda$.
Consequently, the function 
\[
f(j) = \lim_\Lambda f_Z(j)
\]
determines a colouring of $X^{\hookleftarrow F}$ with $m$ colours. \par
Now let $\iota\colon G\hookleftarrow X$ be an arbitrary
isometric embedding, and let $k\in\{1,2,\dots,m\}$ be a colour.
For every $Z\in [\iota(G),\infty)$ choose, 
using ($\ast$), an isometric embedding $j_{Z,k}\colon F\hookrightarrow G$ 
with no element in the 
$\e$-neighbourhood of $\iota\circ j_{Z,k}$, formed in $Z^{\hookleftarrow F}$, 
being of $f_Z$-colour $k$.
For every $x\in F$ define $j_k(x)=\lim_\Lambda j_{Z,k}(x)\in G$.
(The metric space $G$ is finite.) This $j_{k}$ is an isometric embedding
of $F$ into $G$ with the property that the $\e$-neighbourhood of 
$\iota\circ j_k$ formed in all of $X^{\hookleftarrow F}$ 
contains no elements of colour $k$.  
Thus, $\neg$(i) is established. 
\par
(ii) $\Rightarrow$ (i): evident. 
\end{proof}

\begin{corol}
Let $X$ and $Y$ be two metric spaces having, up to isometry, the
same finite metric subspaces. 
If $X$ has property $R$, then so does $Y$. \qed
\label{upto}
\end{corol}

\begin{thm}
Let $X$ and $Y$ be two $\omega$-homogeneous metric spaces,
having, up to isometry, the same finite metric subspaces.
Then the topological group $\Iso(X)$ has the fixed point on compacta property
if and only if the topological group $\Iso(Y)$ does. 
\label{ifandonly}
\end{thm}

\begin{proof} 
Combine Theorem \ref{ramsey} and Corollary \ref{upto}.
\end{proof}

We can finally deduce from Theorem \ref{ifandonly} and Theorem \ref{fpc}
the following result, which is the {\it raison d'\^etre} of the article.

\begin{thm}
Let $U$ be a generalized Urysohn metric space. 
If $U$ is $\omega$-homogeneous,
then the group $\Iso(U)$ has the fixed point on compacta property.
\qed
\label{generurysohn}
\end{thm}

Modulo Uspenskij's Theorem \ref{uspenskij}, the above Theorem
implies the following.

\begin{corol}
Every topological group embeds, as a topological subgroup, into
an extremely amenable topological group, that is, a topological
group with the fixed point on compacta property. \qed
\label{main}
\end{corol}

Even the following appears to be a new result.

\begin{corol}
Every topological group embeds, as a topological subgroup, into
an amenable topological group. \qed
\end{corol}

\subsection{Groups of isometries of discrete Urysohn spaces} 
Here we will demonstrate
how Theorem \ref{ramsey} can be used to show
the absence of the fixed point on compacta property in the case
where the $\omega$-homogeneous metric space in question
fails the `strong' version of Ramsey-type property.

\begin{thm}
The group of isometries of the discrete Urysohn
metric space $\U^{\e \Z}$ does not have the fixed point on compacta property.
\label{doesnot}
\end{thm}

\begin{proof}
Denote by $\{a,b\}$ the two-element metric space with $d(a,b)=\e$.
Partition the set $(\U^{\e \Z})^{\hookleftarrow\{a,b\}}$ of all
isometric embeddings of $\{a,b\}$ into $\U^{\e\Z}$ into two disjoint
subsets $A,B$ in such a way that whenever an injection
$i\colon\{a,b\}\hookrightarrow (\U^{\e\Z})$ is in $A$, the 
`flip' injection $i\circ \sigma_2$ is in $B$, and vice versa.
Since the space $\U^{\e \Z}$ is $\e$-discrete, the $\e$-neighbourhood
of a subset $X$ is $X$ itself, and
`monochromatic up to $\e$' means in this context simply `monochromatic.'
One concludes that, with respect to the colouring $\{A,B\}$, 
no pair of injections of the form $Y=\{i,i\circ\sigma_2\}$ is monochromatic
up to $\e$, and thus the metric
space $(\U^{\e \Z})^{\hookleftarrow\{a,b\}}$,
upon which the group $\Iso(\U^{\e\Z})$ acts transitively and
continuously by isometries, fails the Ramsey--Dviretzky--Milman property.
\end{proof}

\begin{remark}
The same result holds for discrete Urysohn spaces of bounded
diameter, $\U^{\e\Z}_d$.
In particular, letting $\e=1=d$, we obtain a result
proved by the present author in 
\cite{P1}, Th. 6.5: the group of permutations
$S_\infty$ of an infinite set, equipped with the pointwise topology, is not
extremely amenable. (This result seems to answer in
the negative an old question by Furstenberg discussed in \cite{GrM}.) \par
Notice also that the groups of isometries of infinite, 
$\omega$-homogeneous metric spaces need not be extremely amenable.
\par
The countable metric space $\U^\Z_1$, equipped with the
$\{0,1\}$-valued metric, actually satisfies a `weaker' version of the Ramsey
result, namely the one for finite {\it subspaces,} rather than
for their {\it injections,} and this result is the well-known
Finite Ramsey Theorem.
(Cf. Ex. \ref{secmaster}.)
However, as we have just seen, the group fails the `stronger' 
version for embeddings of finite spaces!
The latter circumstance destroys the extreme
amenability of $S_\infty$. 
\par
Finally notice that the topological group $S_\infty$ is amenable,
because it is approximated from within by an increasing
chain of finite groups of permutations whose union is everywhere dense.
\endrem
\end{remark}

\subsection{\label{last}Deducing Ramsey-type theorems for metric spaces}

By force of Theorem \ref{ramsey},
the immediate corollary --- and in fact an equivalent form ---
of the fixed point on compacta property
of the group $\Iso(\U)$ (Theorem \ref{fpc})
is the following Ramsey-type result.

\begin{corol} Let $F$ be a finite metric space, and
let all isometric embeddings of $F$ into $\U$
be coloured using finitely many colours. 
Then for every finite metric space $G$ and every $\e>0$ 
there is an isometric copy $G'\subset \U$ of $G$ such that
all isometric embeddings of $F$ into $\U$ that factor through
$G$ are monochromatic up to $\e$. \qed
\label{ramsey1}
\end{corol}

By restricting ourselves to considering only $\Iso(F)$-invariant
collections of embeddings of $F$ into $\U$, we arrive at the
following.

\begin{corol} Let $F$ be a finite metric space.
Let all subspaces of the Urysohn space $\U$ isometric to $F$ be
coloured using finitely many colours. 
Then for every finite metric space $G$ and every
$\e>0$ there is a
subspace $G'\subseteq\U$ isometric to $G$ whose subspaces
isometric to $F$ are monochromatic up to $\e$. \qed
\label{ramsey2}
\end{corol}

Applications to spherical spaces are probably more interesting.
(Cf. comments in \cite{Mat} at the bottom of p. 460).
The unit sphere of the infinite-dimensional
Hilbert space $\H$ is an $\omega$-homogeneous metric space, and the
orthogonal group of $\H$ 
with the strong operator topology (that is, the
topology of simple convergence on the sphere) is extremely amenable
\cite{GrM}.
As a corollary, we obtain Ramsey-type results for
the Hilbert sphere. 

\begin{corol} Let $F$ be a finite metric subspace of the unit
sphere $\s^\infty$ in an infinite-dimensional Hilbert space.
Let all isometric embeddings of $F$ into $\s^\infty$ be
coloured using finitely many colours. 
Then for every finite metric subspace $G$ of the sphere and
every $\e>0$ there is
an isometric copy $G'\subset\s^\infty$ of $G$ such that
all isometric embeddings of $F$ into 
$G'$ are monochromatic up to $\e$. \qed
\label{ramsey3}
\end{corol}

\begin{corol} Let $F$ be a finite metric subspace of the unit
sphere $\s^\infty$ in an infinite-dimensional Hilbert space.
Let all subspaces of $\s^\infty$ isometric to $F$ be
coloured using finitely many colours. 
Then for every finite subspace $Y$ of the sphere and every
$\e>0$ there is a
subspace $Y'\subseteq\s^\infty$ isometric to $Y$ whose subspaces
isometric to $F$ are monochromatic up to $\e$. \qed
\label{ramsey4}
\end{corol}

To establish similar corollaries for metric subspaces
of the infinite-dimensional
Hilbert space, we need the following result. 
Notice that amenability of the group $\Iso(\H)$ of affine isometries
of a Hilbert space $\H$ was noted in \cite{dlHV}, p. 47.

\begin{thm} 
The group $\Iso(\H)$ of affine isometries of a Hilbert
space $\H$ of infinite dimension is extremely amenable.
\end{thm}

\begin{proof}
The topological group $\Iso(\H)$ is isomorphic to the semidirect
product $\mathrm{O}(\H)\ltimes \H$ of the full orthogonal group
$\mathrm{O}(\H)$ equipped with the strong operator topology and the additive
group of the Hilbert space $\H$ with the usual norm topology, 
formed with respect to the natural action of $\mathrm{O}(\H)$ on $\H$
by rotations.
(Cf. \cite{dlHV}.) Suppose $\Iso(\H)$ acts continuously
on a compact space $K$.
Since the group $\mathrm{O}(\H)$
(identified with a subgroup of $\Iso(\H)$) is extremely
amenable (\cite{GrM}; cf. also Subsection \ref{newp}), 
it has a fixed point $\kappa\in K$.
The mapping $\H\ni x\mapsto x\cdot\kappa\in K$, where $\H$ is viewed
as a closed normal subgroup of
$\Iso(\H)$, is $\Iso(\H)$-equivariant, continuous, and has everywhere
dense image in $K$, and thus $K$ is an equivariant $\Iso(\H)$-compactification 
of the homogeneous factor-space $\H\cong\Iso(\H)/\mathrm{O}(\H)$.\par
Let $\varphi\colon K\to\R^N$ be an arbitrary continuous function,
$N\in\N$.
Its pull-back, $f(x)=:\varphi(x\cdot\kappa)$, to $\H$ is right
uniformly continuous.
(A standard result in abstract topological dynamics.) 
If $\e>0$ is arbitrary, then for some neighbourhood
$V=V[F;\delta]$ of identity in $\Iso(\H)$ one has
$\abs{f(g(0))-f(h(0))}<\e$ whenever $gh^{-1}\in V$.
Without loss in
generality and slightly 
perturbing the points of $F$ if necessary, one can assume
that elements of $F$ are affinely independent.
Let $x,y\in\H$ be two arbitrary elements with the property
$\norm{x-z}=\norm{y-z}$ for each $z\in F$.
Find an isometric copy
of $F$, say $F'$, such that $F'\cup\{0\}$ is isometric to
$F\cup\{x\}$ (or, equivalently, to $F\cup\{y\}$). 
There is an isometry $g$ of $\H$ taking $F'\cup\{0\}$ to $F\cup\{x\}$,
and an isometry $h$ taking $F'\cup\{0\}$ to $F\cup\{y\}$. 
In particular, $gh^{-1}\vert_F=\mathrm{Id}_F\in V$, and consequently
$\abs{f(x)-f(y)}<\e$.
Thus, the function $f$ is $\e$-constant on every affine
sphere of codimension $\abs F$ having the form
$\{x\in\H\colon \norm{x-z}=r_z, z\in F\}\equiv\cap_{z\in F}\s_{r_z}(z)$. 
Another way to say it is that, up to $\e$, the function $f(x)$ only
depends on the collection of distances $\{\norm{x-z}\colon z\in F\}$.
\par
Now let $g_1,\dots,g_n\in\Iso(\H)$ be an arbitrary collection of isometries.
By slightly perturbing them if necessary, one can assume without loss
in generality that all the vectors $z$ and 
$g_i^{-1}(z)$, $z\in F$, $i=1,2,\dots,n$,
are affinely independent. 
Because of infinite-dimensionality of $\H$, 
every element $x$ of some affine subspace of
$\H$ of finite codimension has the property that 
for every $i=1,2,\dots,n$ and each $z\in F$,
one has $\norm{x-g_i^{-1}(z)}=\norm{x-z}$.
Fix any such $x$.
Then the values of
$f$ at the points $x,g_1(x),g_2(x)$, $\dots,g_n(x)$ differ by less than
$\e$.
Now we can apply Theorem \ref{folk} to
conclude that $K$ has a fixed point for $\Iso(\H)$.
\end{proof}

\begin{corol} Let $F$ be a finite metric subspace of the
infinite-dimensional Hilbert space $\H$.
Let all isometric embeddings of $F$ into $\H$ be
coloured using finitely many colours. 
Then for every finite collection $Y$ of such embeddings
and every $\e>0$ there is a collection of embeddings $Y'$
congruent to $Y$ and monochromatic up to $\e$.  \qed
\label{ramsey5}
\end{corol}

\begin{corol} Let $F$ be a finite metric subspace of an infinite-dimensional
Hilbert space $\H$.
If all subspaces of $\H$ isometric to $F$ are coloured
using finitely many colours, then for every finite subspace $G$ of $\H$
and every $\e>0$ there is an isometric copy $G'$ of $G$ in $\H$ such
that all subspaces of $G'$ isometric to $F$ are monochromatic up
to $\e$. \qed
\end{corol}

\section{Concluding remarks}
In this article we have investigated some
relationships inside the following triangle:
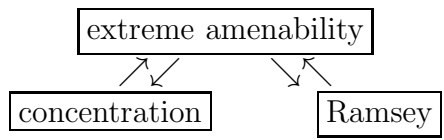
\[
\begin{matrix}
 \mbox{\fbox{extreme amenability}}   \\
 \nearrow\swarrow \phantom{xxxx} \searrow\nwarrow  \\
\mbox{\fbox{concentration}} \phantom{xxxxxx} \mbox{\fbox{Ramsey}} 
\end{matrix}
\]

Deeper explorations of the Ramsey--Milman phenomenon in topological
transformation groups require discovering situations in which a
`phase transition' between concentration and dissipation occurs in
families of topological groups / dynamical systems. (Cf. \cite{BC}.)
It could be, for example, that a solution
to Glasner's problem on the existence of a minimally almost periodic
group topology on the integers without the fixed point on compacta
property \cite{Gl} lies namely in this direction. 

In connection with the Banach--Mazur problem  (cf.  
\cite{Cab}), it could be worth investigating the
fixed point on compacta property for the groups of isometries of
separable Banach spaces admitting a transitive norm.

Finally, we do not know if the results of Section \ref{secappl}
can be put in direct connection with 
the Euclidean Ramsey theory \cite{Graham}.

\newpage
$\,$
\vskip 3cm
\centerline{\bf CORRIGENDUM}
\vskip 1cm            

As kindly pointed out to me by C. Ward Henson, the proof of one of the main
technical results (Theorem 3.2) in the paper is flawed.
To quote from his message:
``Suppose $a,b$ are elements of $F = F^{(3)}_{m+n}$ that have the minimum distance
$\delta$ from each other in the $\rho'$ metric, and let w be any word in 
$F$.  Since
the metric d is bi-invariant, the conjugate
$v = wab^{-1}w^{-1}$ of $ab^{-1}$ has $d$-distance $\delta$ from the identity.
But it seems clear that the reduced length of $v$ could be made arbitrarily
large by choosing w correctly.  This contradicts what you claim in (8).''

Fortunately, the result is not particularly deep, and here is a corrected
proof of the statement.

As in \cite{urysohn}, we 
say that a metric space $X$ is {\bf indexed} by a set $I$ if
there is a surjection $f_X\colon I\to X$.
We will call the pair
$(X,f_X)$ an {\bf indexed metric space.} Two metric
spaces, $X$ and $Y$, indexed with the same set $I$ are 
$\e${\bf -isometric} if for every $i,j\in I$ the distances
$d_X(f_X(i),f_X(j))$ and $d_Y(f_Y(i),f_Y(j))$ differ by at most $\e$.

Here is the result in question.
\\[5mm]
\noindent{\bf Theorem 3.2.}
{\it
Let $g_1,\dots,g_m$ be a finite family of isometries of a
metric space $X$.
Then for every $\e>0$ and every finite collection
$x_1,\dots, x_n$ of elements of $X$ there exist a finite metric space
$\widetilde X$, elements $\tilde x_1,\dots,\tilde x_n$ of $\tilde X$,
and isometries $\tilde g_1,\dots,\tilde g_m$ of $\tilde X$ such that
the indexed metric spaces 
$\{g_j\cdot x_i\colon  i=1,2,\dots,n, j=1,2,\dots,m\}$ 
and $\{\tilde g_j\cdot \tilde x_i\colon  i=1,2,\dots,n, j=1,2,\dots,m\}$
are $\e$-isometric.}

\begin{proof} 
Without loss in generality, we can assume that $X$ is separable
(in fact, even countable). 
Such an $X$ can be $g$-embedded into the Urysohn space $\U$
(see \cite{usp}, or else Prop. 4.1 in \cite{urysohn}),
and therefore we can further assume that $X=\U$. 

Choose any element $\xi\in\U$ and isometries $g_{m+1},g_{m+2},\ldots,g_{m+n}$
of $\U$
with the property $g_{m+i}(\xi) = x_i$, $i=1,2,\ldots,n$.

Denote by $F_{m+n}$ the free non-abelian group on generators
$g_1,\ldots,g_m$, $g_{m+1},\ldots,g_{m+n}$. The group 
$F_{m+n}$ acts on $\U$ by isometries.

The formula
\[d(g,h):= d_\U(g(\xi),h(\xi)),~~g,h\in F_{m+n},\]
where $d_\U$ denotes the metric on the Urysohn space, defines
a left-invariant pseudometric $d$ on the group $F_{m+n}$:
\begin{eqnarray*}
d(xg,xh) &=& d_\U(xg(\xi),xh(\xi))\\
&=& d_\U(x(g(\xi)),x(h(\xi)))\\
&=& d_\U(g(\xi),h(\xi))\\ &=& d(g,h).
\end{eqnarray*}
The indexed
metric subspace $\{g_j\cdot x_i\colon  i=1,2,\dots,n, j=1,2,\dots,m\}$ 
of $\U$ is
isometric to the metric subspace 
$\{g_j\cdot g_{m+i}\colon  i=1,2,\dots,n, j=1,2,\dots,m\}$ of 
$(F_{m+n},d)$. Indeed,
\begin{eqnarray*}
d_\U(g_j\cdot x_i,g_k\cdot x_l) &=& d_\U(g_j\cdot g_{m+i}(\xi),
g_k\cdot g_{m+l}(\xi))
\\
&=& d(g_j g_{m+j},g_k g_{m+l}).
\end{eqnarray*}

Notice also that the latter subspace is contained
in the set $F_{m+n}^{(2)}$ of all words having reduced length $\leq 2$.
(The reduced length will always mean that with regard to the
generators $g_1,\ldots,g_m$, $g_{m+1},\ldots,g_{m+n}$.)

By adding to $d$ a left-invariant metric on $F_{m+n}$ taking sufficiently
small values on pairs of elements of $F_{m+n}^{(2)}$, we can assume
without loss in generality that $d$ is a left-invariant metric on $F_{m+n}$.
(For instance, add a metric whose only non-zero value is $\e/3$.)

Form the Cayley graph $\Gamma$
of $F_{m+n}$ with
regard to the set of generators $Y = F_{m+n}^{(4)}$. Vertices of
$\Gamma$ are elements of $F_{m+n}$, and $x,y\in F_{m+n}$ are adjacent
if and only if $x^{-1}y\in F_{m+n}^{(4)}$. This graph is connected.
Now make $\Gamma$ into a weighted graph by assigning to an edge 
$(a,b)$, $a^{-1}b\in F_{m+n}^{(4)}$ the value $d(a,b)\equiv d(a^{-1}b,e)$.

Denote by $\rho$ the path metric of the weighted graph $\Gamma$. Its
value for $x,y\in F_{m+n}$ is given by 
\begin{equation}
\rho(x,y) = \inf\sum_{i=0}^{N-1} d(a_i,a_{i+1}),
\label{minimum}
\end{equation}
where the infimum is taken over all natural $N$ and all finite sequences
$x=a_0, a_1,\ldots,a_{N-1},a_N=b$, with the property 
$a_i^{-1}a_{i+1}\in F_{m+n}^{(4)}$ for all $i$. 

It is easily seen that $\rho$ is a left-invariant metric on the group 
$F_{m+n}$.

Generally, $\rho\geq d$, but restrictions of $\rho$ and $d$ to
$F_{m+n}^{(2)}$ coincide.

If one denotes by $\delta>0$ the minimal value of $d(a,b)$ as
$a,b\in F_{m+n}^{(2)}$ and $a\neq b$,
then $\rho(x,y)\geq \delta d_w(x,y)$, where $d_w$ denotes the word metric
with respect to the set of generators $Y = F_{m+n}^{(4)}$. 
In particular, if an $x\in F_{m+n}$ has reduced length $l=l(x)$, then 
$d_w(x)\geq l/4$ and accordingly $\rho(x,e)\geq \delta l/4$.
(As a consequence, the infimum in Eq. (\ref{minimum}) is always achieved.)

Let $\Delta$ denote the maximal value of the metric $\rho$ between
pairs of elements of $F_{m+n}^{(2)}$.
Choose a natural number $N$ so large that $\delta(N-4)/4\geq\Delta$, for
instance, set $N=4\lfloor \Delta/\delta\rfloor + 4$.

Every free group is residually finite, that is, admits a separating
family of homomorphisms into finite groups.
(Cf. e.g. \cite{Hall+}, Ch. 7, exercise 5.)
Using this fact, choose a normal
subgroup $H\triangleleft F_{m+n}$ of finite index so that 
$H\cap F_{m+n}^{(N)} =\{e\}$.

The formula
\begin{eqnarray}
\label{infimum}
\bar \rho(xH,yH) &:=& \inf_{h_1,h_2\in H} \rho(xh_1,yh_2) \\
&\equiv& \inf_{h_1,h_2\in H} \rho(h_1x,h_2y) \nonumber\\
&\equiv& \inf_{h\in H} \rho(hx,y)  \nonumber
\end{eqnarray}
defines a left-invariant pseudometric on the finite factor-group
$F_{m+n}/H$. The triangle inequality follows from the fact that,
for all $h'\in H$,
\begin{eqnarray*}
\bar \rho(xH,yH) &=&
\inf_{h\in H} \rho(hx,y) \\
&\leq &
\inf_{h\in H} [ \rho(hx,h'z) +  \rho(h'z,y) ]\\
&=& \inf_{h\in H} \rho(hx,h'z) +  \rho(h'z,y) \\
&=& \inf_{h\in H} \rho(h'^{-1}hx,z) +  \rho(h'z,y) \\
&=& \bar\rho(xH,zH) + \rho(h'z,y),
\end{eqnarray*}
and the infimum of the r.h.s. taken over all  $h'\in H$ equals
$\bar\rho(xH,zH) + \bar\rho(zH,yH)$. Left-invariance of $\bar\rho$ is obvious.
%\begin{eqnarray*}
%\bar \rho(zxH,zyH) &=& \inf_{h_1,h_2\in H} \rho(zxh_1,zyh_2) \\
%&=& \inf_{h_1,h_2\in H} \rho(xh_1,yh_2) \\
%&=& \bar \rho(xH,yH).
%\end{eqnarray*}

Let $x,y\in F_{m+n}^{(2)}$. Closely approximate the infimum 
in Eq. (\ref{infimum}) by some value $\rho(xh_1,yh_2)$ with 
$h_1,h_2\in H$, 
then 
\[\rho(xh_1,yh_2) =\rho(y^{-1}x h_1 x^{-1}y\cdot y^{-1}x, h_2)
= \rho(y^{-1}x, h_3),
\]
where $h_3 = y^{-1}x h_1^{-1} x^{-1}y h_2\in H$. 

The value $\rho(y^{-1}x, h_3)$, $h_3\in H$,
cannot get smaller than $d(y^{-1}x,e) = d(x,y)$. Indeed, 
unless $h_3= e$ (in which case $\rho(y^{-1}x, h_3) = \rho(x,y) =
d(x,y)$), one has
$l(h_3)\geq N$ and so the word distance from $y^{-1}x$ to $h_3$ is at least
$N-4$, and $\rho(y^{-1}x,h)\geq \delta(N-4)/4\geq\Delta\geq d(x,y)$.

We conclude: the restriction of the factor-homomorphism
\[\pi\colon F_{m+n}\ni x\mapsto xH \in F_{m+n}/H\]
to $F_{m+n}^{(2)}$ is an isometry. 

One can now perturb the pseudometric on $F_{m+n}/H$ by adding to it
a left-invariant
metric taking very small values (e.g. taking the only 
non-zero value $\e/3$) so as to replace $\bar\rho$ 
with a left-invariant metric, $\tilde\rho$. 

Take now $\tilde X =
(F_{m+n}/H,\tilde\rho)$, $\tilde x_i =\pi(g_{m+i})\in\tilde X$, $i=1,2,\ldots,n$,
and let $\tilde g_j$ be left translates made by the elements
$\pi(g_j)$, $j=1,2,\ldots,m$, in the finite group $F_{m+n}/H$.
The indexed
metric space  
$\{g_j\cdot g_{m+i}\colon  i=1,2,\dots,n, j=1,2,\dots,m\}$ of 
$(F_{m+n},d)$ is $\e$-isometric to the metric subspace 
$\{\pi(g_j)\cdot \pi(g_{m+i})\colon  i=1,2,\dots,n, j=1,2,\dots,m\}$ of 
$(F_{m+n}/H,\tilde\rho)$. Consequently,
 the conclusion of the Theorem is verified.
\end{proof}

\noindent
{\bf Remark.}
Prof. Henson has also pointed out to me that in the particular case of
path metric spaces associated to a
graph the above result (Theorem 3.2)
follows from earlier results by Hrushovski \cite{hrushovski}.
\\[4mm]
{\bf Acknowledgements.}
The author is most grateful to Prof. C. Ward Henson for his 
comments, as well as 
to all other participants of the Research Among Peers (RAP) seminar series
``Polish group actions and extremely amenable groups'' otganized 
by Prof. Slawomir Solecki ind the Department of Mathematics,
University of Illinois at Urbana-Champaign 
in Spring 2003.

\end{document}